\documentclass[11pt]{amsart}
\usepackage{graphicx}
\usepackage{amscd}

\usepackage{fixltx2e}
\usepackage[utf8]{inputenc}
\usepackage{amssymb, latexsym,amsthm,amsfonts, amsbsy, amsmath, mathrsfs}
\usepackage{bm}                          
\usepackage{verbatim}                    
\usepackage[all, knot]{xy}               
        \xyoption{arc}
        \xyoption{web}                   
\makeatletter                            
\def\MT@register@subst@font{\MT@exp@one@n\MT@in@clist\font@name\MT@font@list
 \ifMT@inlist@\else\xdef\MT@font@list{\MT@font@list\font@name,}\fi}
\makeatother

\newtheorem{theorem}{Theorem}
\newtheorem{lemma}[theorem]{Lemma}

\newtheorem{corollary}[theorem]{Corollary}

\newtheoremstyle{inproof}
{3pt}
{3pt}
{\normalfont}
{}
{\bfseries}
{:}
{.5em}
{}

\theoremstyle{inproof}
\newtheorem{claim}{Claim}

\theoremstyle{definition}
\newtheorem{definition}[theorem]{Definition}

\newtheorem{remark}[theorem]{Remark}

\newcommand{\rig}{\rightarrow}
\newcommand{\mrig}{\mathrel{-\!\!\!\!\!\rightarrow}}

\newcommand{\Rig}{\Rightarrow}

\newcommand{\bcdw}{\mathbin{\boldsymbol\cdot}}
\newcommand{\bcdn}{\mbox{\boldmath{$\cdot$}}}

\newcommand{\seteq}{\mathrel{\mbox{\,:\!}=\nolinebreak }\,}

\newcommand{\up}{\mathop{\,\uparrow}}
\newcommand{\down}{\mathop{\,\downarrow}}

\newcommand{\sbA}{{\boldsymbol{A}}}
\newcommand{\sbB}{{\boldsymbol{B}}}
\newcommand{\sbC}{{\boldsymbol{C}}}
\newcommand{\sbD}{{\boldsymbol{D}}}
\newcommand{\sbE}{{\boldsymbol{E}}}

\newcommand{\sbJ}{{\boldsymbol{J}}}

\newcommand{\sbS}{{\boldsymbol{S}}}

\newcommand{\sbX}{{\boldsymbol{X}}}

\newcommand{\ld}{{\backslash}}

\newcommand{\ov}{\overline}

\newcommand{\monus}{\mbox{$\frac{\:{\bf \displaystyle .}\:}{\mbox{}}$}}

\let\class=\mathsf                                 
\bmdefine{\A}{A}                                   
\bmdefine{\B}{B}
\bmdefine{\C}{C}
\bmdefine{\D}{D}
\bmdefine{\Fm}{Fm}                                 

\bmdefine{\leibniz}{\varOmega}                       

\DeclareMathOperator{\Fg}{\mathit{Fg}}

\DeclareMathOperator{\Cpu}{Cpu}
\DeclareMathOperator{\Sg}{\mathit{Sg}}

\providecommand{\monus}{
	\mathbin{
		\vphantom{+}
		\text{
			\mathsurround=0pt 
			\ooalign{
				\noalign{\kern-.35ex}
				\hidewidth$\smash{\cdot}$\hidewidth\cr 
				\noalign{\kern.35ex}
				$-$\cr 
			}%
		}%
	}%
}
\newcommand{\alg}[1]{{\boldsymbol{#1}}}
\newcommand{\spc}[1]{{\boldsymbol{#1}}}

\newcommand{\cl}[1]{{\mathsf{#1}}}
\newcommand{\id}{{\textup{id}}}

\newcommand{\bdot}{\mathbin{\boldsymbol{\cdot}}}

\DeclareMathOperator{\depth}{depth}

\newcommand{\upclose}{{\uparrow}}
\newcommand{\downclose}{{\downarrow}}

\makeatletter
\providecommand*{\Dashv}{%
	\mathrel{%
		\mathpalette\@Dashv\vDash
	}%
}
\newcommand*{\@Dashv}[2]{%
	\reflectbox{$\m@th#1#2$}%
}
\makeatother
\makeatletter

\title[Epimorphisms and Subidempotent Residuated Structures]{Epimorphisms in Varieties of Subidempotent Residuated Structures}
\author{T.\ Moraschini}
\address{Institute of Computer Science, Academy of Sciences of the Czech Republic, Pod Vod\'{a}renskou v\v{e}\v{z}\'{i} 2, 182 07 Prague 8, Czech Republic.}
\email{moraschini@cs.cas.cz}
\author{J.G.\ Raftery}
\address{Department of Mathematics and Applied Mathematics,
 University of Pretoria,
 Private Bag X20, Hatfield,
 Pretoria 0028, South Africa}
\email{{james.raftery@up.ac.za} \quad {jamie.wannenburg@up.ac.za}}
\author{J.J.\ Wannenburg}
\keywords{Epimorphism, residuated lattice, Brouwerian algebra, Heyting algebra, De Morgan monoid, Esakia space, substructural logic,
relevance logic, Beth definability.}
\subjclass[2010]{{}\\ \indent Primary: 03B47, 03G25, 06F05.  Secondary: 03B55, 06D20.}
\thanks{This work received funding from the European Union's Horizon 2020 research and innovation programme under the Marie Sklodowska-Curie grant
agreement No~689176 (project ``Syntax Meets Semantics: Methods, Interactions, and Connections in Substructural logics").  The first author was also supported by the project CZ.02.2.69/0.0/0.0/17\_050/0008361, OPVVV M\v{S}MT, MSCA-IF Lidsk\'{e} zdroje v teoretick\'{e} informatice.
The second author was supported in part by the National Research Foundation of South Africa (UID 85407).
The third author was supported by the DST-NRF Centre of Excellence in Mathematical and Statistical Sciences (CoE-MaSS), South Africa.  Opinions expressed and conclusions arrived at are those of the authors and are not necessarily to be attributed to the CoE-MaSS}

\begin{document}

\begin{abstract}
A commutative residuated lattice $\sbA$ is said to be \emph{subidempotent} if the lower bounds of its neutral element $e$ are idempotent (in which case they naturally constitute a Brouwerian algebra $\sbA^-$).
It is proved here that epimorphisms are surjective in
a variety $\cl{K}$ of such algebras $\sbA$
(with or without
involution), provided that each finitely subdirectly irreducible algebra $\sbB\in\cl{K}$ has
two properties: (1)~$\sbB$ is generated by
lower bounds of $e$, and (2)~the poset of prime
filters of $\sbB^-$ has finite depth.  Neither (1) nor (2) may be dropped.
The proof adapts to the presence of bounds.  The result generalizes some
recent findings of G.\ Bezhanishvili and the first two authors concerning epimorphisms in varieties
of Brouwerian algebras, Heyting algebras and Sugihara monoids, but its scope also encompasses a range of
interesting varieties of De Morgan monoids.
\end{abstract}

$\mbox{}$\\{\vspace{-8mm}}

\maketitle

\makeatletter
\renewcommand{\labelenumi}{\text{(\theenumi)}}
\renewcommand{\theenumi}{\roman{enumi}}
\renewcommand{\theenumii}{\roman{enumii}}
\renewcommand{\labelenumii}{\text{(\theenumii)}}
\renewcommand{\p@enumii}{\theenumi(\theenumii)}
\makeatother

\allowdisplaybreaks

\section{Introduction}

In a variety of algebras, if a homomorphism is surjective, then it is an epimorphism, but the converse need
not hold.  Indeed, rings and distributive lattices each form varieties in which non-surjective
epimorphisms arise.  As it happens, this reflects the absence of unary terms defining multiplicative
inverses in rings, and complements in distributive lattices, despite the uniqueness of those entities
when they exist.

Such constructs are said to be \emph{implicitly} (and not \emph{explicitly}) \emph{definable}.
In a variety of logic, they embody implicitly definable propositional functions that cannot be
explicated in the corresponding logical syntax, and Beth-style `definability
properties' \emph{preclude} phenomena of this kind.  (The allusion is to E.W.\ Beth's theorems for
classical propositional and predicate logic in \cite{Bet53}.)

In particular, when a logic $\mathbf{L}$ is
algebraized, in the sense of \cite{BP89}, by a variety $\cl{K}$ of algebras, then the
\emph{ES property} for $\cl{K}$---i.e., the demand that \emph{all} epimorphisms in $\cl{K}$ be surjective---amounts
to the so-called \emph{infinite Beth definability property} for $\mathbf{L}$.  The most general version of this `bridge theorem'
was formulated and proved by Blok and Hoogland \cite[Theorems~3.12,\,3.17]{BH06} (also see \cite[Theorem~7.6]{MRW3} and
the antecedents cited in both papers).

In this situation, the subvarieties of $\cl{K}$ algebraize the axiomatic extensions of $\mathbf{L}$, but the
ES property need not persist in subvarieties.  It is therefore a well-motivated (but often nontrivial) task to determine
which subvarieties of $\cl{K}$ have surjective epimorphisms.  The present paper addresses this question in the
context of residuated structures, i.e., the algebraic models of substructural logics (see \cite{GJKO07,KO10}).

We consider algebras called \emph{SRLs}, which are residuated lattice-ordered commutative monoids
satisfying the \emph{subidempotent law} $x\leqslant e\Longrightarrow x\bcdw x=x$ (where $\bcdw$ is the monoid operation
and $e$ its neutral element)
and possibly
equipped with an involution (\emph{SIRLs}) and/or lattice-bounds.  They need not be
\emph{integral}, i.e.,
$e$ need not be the greatest element.  They include De Morgan
monoids, i.e., the algebraic models of the relevance logic $\mathbf{R^t}$ of \cite{AB75}.
The \emph{negative cone} of an S[I]RL $\sbA$, which comprises the lower bounds of $e$, may be given
the structure of a Brouwerian or Heyting algebra $\sbA^-$, to which the Esakia duality of
\cite{Esa85} applies.  In particular, the \emph{depth} of $\sbA$ may be defined as that of $\sbA^-$.

Heyting and Brouwerian algebras model intuitionistic propositional logic and its positive fragment.
A result of Kreisel \cite{Kre60} shows (in effect) that every variety consisting of such algebras
has a \emph{weak} form of the ES property, whereas Maksimova established that only finitely many
enjoy a
certain \emph{strong} form;
see \cite{GM05,Mak00,Mak03}.  Uncountably
many varieties of Heyting [Brouwerian] algebras have finite depth \cite{Kuz75}, and all of these
have surjective epimorphisms; the latter claim was proved recently by G.\ Bezhanishvili and the
first two authors \cite{BMR17}, using Esakia duality.
On the other hand, the ES property fails in uncountably many further varieties
of Heyting [Brouwerian] algebras \cite{MW}.

Relational duals for non-integral S[I]RLs are sometimes available \cite{Urq96}, but they are rather complicated.
Also, the functor that constructs negative cones (and restricts morphisms accordingly) is not a category equivalence,
except in quite special cases; see \cite{FG19,GP,GR12,GR15}.  Partly for these reasons, we cannot systematically
reduce ES problems for arbitrary varieties of
S[I]RLs to an examination
of negative cones.

Nevertheless, $\sbA^-$ contains enough
information about $\sbA$ to facilitate the main result of this paper (Theorem~\ref{thm:SIRLhasES}), which
is a sufficient condition for the surjectivity of epimorphisms.  It states that,
in a variety
$\cl{K}$
of S[I]RLs, epimorphisms will be surjective if each finitely subdirectly irreducible member of $\cl{K}$
has finite depth and is generated by elements of its negative cone.
Neither hypothesis may be dropped.

The assumptions of
Theorem~\ref{thm:SIRLhasES} persist in subvarieties and under varietal joins, so the result is
labour-saving.  Apart from generalizing the aforementioned findings of \cite{BMR17}, it
settles the question of epimorphism-surjectivity for many interesting varieties of De Morgan monoids
(see Sections~\ref{reflections} and \ref{section:applications}),
yielding definability results for a corresponding
range
of relevance logics.

\medskip

\noindent{\bf Notation.}
In an indicated partially ordered set $\langle X;\leqslant\rangle$,
we define
\[
\textup{$\up{x}=\{y\in X\colon x\leqslant y\}$
\ and \ $\up U=\mbox{\scriptsize $\bigcup$}_{u\in U}\up{u}$,}
\]
for $U\cup\{x\}\subseteq X$, and if $U=\up U$, we call $U$ an {\em up-set\/} of $\langle X;\leqslant\rangle$.
We define $\down x$ and $\down U$
dually.  For $x,y\in X$, the notation $x\prec y$ (`$y$ covers $x$') signifies that $x<y$ and there is no $z\in X$ such that $x<z<y$.

If $Y\subseteq X$ and $x\in X$, we sometimes need to refer to $Y\cap\up{x}$, which we then denote as $\upclose^Y x$, even if $x\notin Y$.

The universe of an algebra $\sbA$ is denoted by $A$.  Thus, the congruence lattice $\boldsymbol{\mathit{Con}}\,\sbA$
of $\sbA$ has universe $\mathit{Con}\,\sbA$.
For $\emptyset\neq X\subseteq A$,
the subalgebra of $\sbA$ generated by $X$ is
denoted by $\boldsymbol{\mathit{Sg}}^\sbA X$ (and its universe by
$\Sg^\sbA X$).

The class operator symbols
$\mathbb{H}$, $\mathbb{S}$, $\mathbb{P}$
and $\mathbb{P}_\mathbb{U}$
stand, respectively, for closure under
homomorphic images, subalgebras,
direct
products
and ultra\-products,
while $\mathbb{V}$
denotes varietal
generation, i.e., $\mathbb{V}=\mathbb{HSP}$.
We abbreviate $\mathbb{V}(\{\sbA\})$ as $\mathbb{V}(\sbA)$.

Recall that an algebra $\sbA$ is \emph{finitely subdirectly irreducible} (briefly, \emph{FSI}) iff
its identity relation $\textup{id}_A$ is meet-irreducible
in
${\boldsymbol{\mathit{Con}}}\,\sbA$.
For a variety $\cl{K}$, we use $\cl{K}_{FSI}$ to denote the class of all FSI members of $\cl{K}$.
Of course, $\cl{K}=\mathbb{V}(\cl{K}_{FSI})$, by the Subdirect Decomposition Theorem.

\section{Epimorphisms}

Given a class $\cl{K}$ of similar algebras, a \emph{$\cl{K}$-morphism} is a homomorphism $\textup{$f \colon \alg{A} \mrig \alg{B}$}$,
where $\alg{A},\alg{B} \in \cl{K}$.  It is called a
\emph{$\cl{K}$-epimorphism} provided that, whenever
$g,h\colon\alg{B}\mrig\alg{C}$ are $\cl{K}$-morphisms
with $g \circ f = h \circ f$, then $g = h$.
Clearly, surjective $\cl{K}$-morphisms are $\cl{K}$-epimorphisms.  We say that
$\cl{K}$ has the \emph{epimorphism-surjectivity (ES) property} if all $\cl{K}$-epimorphisms are surjective.

A subalgebra $\alg{D}$ of an algebra $\alg{E}\in\cl{K}$ is said to be \emph{$\cl{K}$-epic} (in $\alg{E}$)
if every $\cl{K}$-morphism with
domain $\alg{E}$ is determined by its restriction to $\alg{D}$.  (This means that the inclusion map $\alg{D}\mrig\alg{E}$
is a $\cl{K}$-epimorphism, assuming that $\alg{D}\in\cl{K}$.)
Thus,
a $\cl{K}$-morphism is a $\cl{K}$-epimorphism iff its image is a $\cl{K}$-epic subalgebra of its co-domain.
And, when $\cl{K}$ is closed under subalgebras,
it has the ES property iff none of its members has a $\cl{K}$-epic proper subalgebra.

Recall that if $\sbA$ is a subalgebra of an algebra $\sbB$, and $\mu\in \mathit{Con}\,\sbB$,
then the relation $\mu|_A\seteq A^2\cap\mu$ is a congruence of $\sbA$.

\begin{lemma}
\label{lem:epi}
Let\/ $\cl{K}$ be a variety of algebras, and\/
$\sbA$
a\/ $\cl{K}$-epic subalgebra of\/ $\sbB\in\cl{K}$\textup{.}
Then, for any\/ $\mu\in\mathit{Con}\,\sbB$\textup{,}
the map\/ $f \colon \alg{A}/(\mu|_A) \mrig \alg{B}/\mu$ defined by\/ $a/(\mu|_A) \mapsto a/\mu$ is an injective\/ $\cl{K}$-epimorphism.
\end{lemma}

\begin{proof}
As $\cl{K}$ is a variety, $\alg{B}/\mu,\,\alg{A}/(\mu|_A)\in\cl{K}$.
Let $i\colon\sbA\mrig\sbB$ be the inclusion homomorphism and
$q\colon \alg{B} \mrig \alg{B}/\mu$ the surjective homomorphism $b\mapsto b/\mu$.
Clearly,
$\mu|_A$ is the kernel of the homomorphism $q \circ i \colon \alg{A} \mrig \alg{B}/\mu$, so $f$ is an injective $\cl{K}$-morphism.

Suppose
$g,h \colon \alg{B}/\mu \mrig \alg{C}\in\cl{K}$ are homomorphisms with
$g \circ f = h \circ f$.  Then
$g \circ q$ and $h \circ q$ are homomorphisms from $\alg{B}$ to $\alg{C}$.  For each $a \in A$,
$$ g(q(a)) = g(a/\mu) = g(f(a/(\mu|_A))) = h(f(a/(\mu|_A))) = h(q(a)),$$
i.e., $g \circ q \circ i = h \circ q \circ i$.  Therefore, $g \circ q = h \circ q$, as $i$ is a $\cl{K}$-epimorphism.
Since $q$ is surjective, it follows that $g=h$, as required.
\end{proof}

A variety $\cl{K}$ is said to have \emph{EDPM} if it is congruence distributive and $\cl{K}_{FSI}$ is a universal class (i.e.,
subalgebras and ultraproducts of FSI members of $\cl{K}$ are FSI).  The acronym stands for `equationally definable principal meets'
and is motivated by other characterizations of the notion in \cite{BP86,CD90}.

\begin{theorem}
\label{thm:ESWitnessedByFSI}
\textup{(Campercholi \cite[Theorem~6.8]{Cam18})}
If a congruence permutable variety\/ $\cl{K}$ with EDPM lacks the ES property, then some FSI member of\/
$\cl{K}$
has a\/ $\cl{K}$-epic proper subalgebra.
\end{theorem}

\section{Residuated Structures}\label{residuated structures section}

General information about residuated structures (and their connection with substructural logics) can be
found in \cite{GJKO07}.
Here, we recall the basic definitions and facts that will be relied on below.
\begin{definition}
An algebra $\sbA=\langle A;\wedge,\vee,\bcdw,\rig,e\rangle$ will be called a (commutative)
\emph{subidempotent residuated lattice},
or briefly an {\em SRL}, if $\langle A;\wedge,\vee\rangle$ is a lattice,
$\langle A;\bcdw,e\rangle$ is a commutative monoid and $\rig$ is a binary operation (called the \emph{residual})
such that $\sbA$ satisfies
\begin{align}
& x\bcdw y\leqslant z\;\Longleftrightarrow\;x\leqslant y\rig z \quad \text{(the \emph{law of residuation})} \label{residuation}\\
& x\leqslant e\;\Longrightarrow\;x=x^2 \textup{ \,$(\!\seteq x\bcdw x)$}
\quad \text{(\emph{subidempotence})}, \label{subidempotence}
\end{align}
where $\leqslant$ is the lattice order.  An equational paraphrase of (\ref{subidempotence}) is
\[
(x\wedge e)^2=x\wedge e.
\]
\end{definition}
Every SRL
satisfies the postulates below.
Here and subsequently,
${x\leftrightarrow y}$ abbreviates ${(x\rig y)\wedge (y\rig x)}$.
\begin{align}
& x\bcdw(y\vee z)=(x\bcdw y)\vee(x\bcdw z)\label{fusion distributivity}\\
& x\rig(y\wedge z)=(x\rig y)\wedge(x\rig z)
\label{and distributivity}\\
&
(x\vee y)\rig z=(x\rig z)\wedge(y\rig z)\label{or distributivity}\\
& x\leqslant y \;\Longleftrightarrow\; e\leqslant x\rig y \label{t order}\\
& x=y \;\Longleftrightarrow\; e\leqslant x\leftrightarrow y \label{t reg}\\
& e\leqslant x\rig x \textup{ \ and \ } e\rig x=x \label{t laws}\\
& (x\leqslant e\;\,\&\,\;y\leqslant e)\;\Longrightarrow\;x\wedge y=x\bcdw y.\label{square increasing cor2}
\end{align}
Of these properties, only
(\ref{square increasing cor2}) relies on subidempotence.
It follows from (\ref{fusion distributivity}) that $\bcdn$ is isotone in both arguments,
and from (\ref{and distributivity}) and (\ref{or distributivity}) that $\rig$ is isotone in its second argument and antitone
in its first.  In particular, the implication $x\leqslant e\;\Longrightarrow\;x^2\leqslant x$ does not rely on (\ref{subidempotence}),
so we could express
(\ref{subidempotence}) as
\[
x\leqslant e\;\Longrightarrow\;x\leqslant x^2.
\]
The key step in the proof of (\ref{square increasing cor2}) is
$x\wedge y\wedge e=(x\wedge y\wedge e)^2\leqslant x\bcdw y$.

An SRL $\sbA$ is said to be
\begin{itemize}
\item \emph{square-increasing} if it satisfies $x\leqslant x^2$;
\item \emph{idempotent} if it satisfies $x=x^2$;
\item \emph{distributive} if its lattice reduct $\langle A;\wedge,\vee\rangle$ is
distributive; and
\item \emph{integral} if $e$ is its greatest element.
\end{itemize}
By (\ref{square increasing cor2}), an SRL
is integral iff its operations $\bcdn$ and $\wedge$
coincide.
\begin{definition}\label{def:brouwerian}
A \emph{Brouwerian algebra} is an integral SRL $\sbA$; it is normally
identified with its reduct $\langle A;\wedge,\vee,\rig,e\rangle$ (in view of the previous remark).
\end{definition}
\noindent
It follows that
a Brouwerian algebra is idempotent, distributive (by (\ref{fusion distributivity})) and determined by its lattice reduct, and that
it satisfies $x\rig e=e=x\rig x$.

\begin{definition}
An \emph{involutive SRL}, briefly an \emph{SIRL}, is the expansion $\sbA$ of an SRL by an \emph{involution}, i.e., a
unary operation $\neg$ such that $\sbA$ satisfies $\neg\neg x=x$ and
$x\rig\neg y=y\rig\neg x$.  \,In these algebras,
we define $f=\neg e$.
\end{definition}
\noindent An SIRL
satisfies
$x\rig y=\neg(x\bcdw\neg y)$
and $\neg x=x\rig f$.
Every SRL can be embedded into an SIRL, in such a way that distributivity and the square-increasing law $x\leqslant x^2$ are preserved
(see \cite{Mey73,GR04} and Section~\ref{reflections} below).

The class of all S[I]RLs is a finitely axiomatizable variety.  It is arithmetical, i.e., congruence distributive and
congruence permutable. (See \cite[p.\,94]{GJKO07}, where a termwise-equivalent formulation is used.)
The variety
of Brouwerian algebras has the ES property \cite{Mak03} (also see \cite{EG81,GM05,Mak00}).

Distributive square-increasing SIRLs---a.k.a.\ \emph{De Morgan monoids}---alge\-braize
the principal relevance logic $\mathbf{R^t}$ of \cite{AB75} (also see \cite{Dun66,MDL74,MRW,MRW2}), while distributive
square-increasing SRLs---a.k.a.\ \emph{Dunn monoids}---algebraize
the `positive' (negationless) fragment of $\mathbf{R^t}$.  Their non-distributive counterparts similarly match
the subsystem
$\mathbf{LR^t}$ of $\mathbf{R^t}$ and its positive fragment.  ($\mathbf{LR^t}$ adds the double-negation axiom to the system
$\mathbf{FL_{ec}}$ of \cite{GJKO07}; also see \cite{Mey01,TMM88}.)
The square-increasing law embodies the logical \emph{contraction} axiom $\textup{$(p\rig(p\rig q))\rig(p\rig q)$}$.
The idempotent De Morgan monoids are called \emph{Sugihara monoids}.  They model the logic $\mathbf{RM^t}$
(a.k.a.\ $\mathbf{R}$-\emph{mingle}), and their structure is better understood than that of other De Morgan monoids;
see \cite{AB75,Dun70,GR12,GR15,OR07}.
In all of these cases, there is a transparent lattice anti-isomorphism between the axiomatic
extensions of the logic and the subvarieties of the model class.

A {\em bounded S[I]RL\/} is the expansion of an S[I]RL by a distinguished element $\bot$, which is the least element of the order, whence $\bot\rig\bot$
is the greatest element.
A \emph{Heyting algebra} is a bounded Brouwerian algebra.
For simplicity, we
defer further discussion of bounded S[I]RLs until Section~\ref{bounds}.

\section{Filters}\label{section:filters}

Let $\sbA$ be an
S[I]RL.
By a \emph{filter} of
$\sbA$, we mean a filter of the lattice $\langle A;\wedge,\vee\rangle$, i.e.,
a non-empty subset $F$ of $A$ that is upward closed
and closed under the binary operation $\wedge$.
It is called a \emph{deductive filter} of $\sbA$ if, moreover, $e\in F$.  In that case, $F$ is a submonoid of
$\langle A;\bcdw,e\rangle$ and, whenever $b\in A$ and $a,a\rig b\in F$, then $b\in F$.  The respective explanations
are that $\sbA$ satisfies
$x\bcdw y\geqslant (x\wedge e)\bcdw (y\wedge e)=x\wedge y\wedge e$ (by (\ref{square increasing cor2}))
and
$x\bcdw(x\rig y)\leqslant y$ (by (\ref{residuation})).

The lattice of deductive filters of $\sbA$ is isomorphic to
${\boldsymbol{\mathit{Con}}}\,\sbA$.
The isomorphism and its inverse are given by
\begin{eqnarray*}
& F\,\mapsto\,\leibniz^\sbA F\seteq\{\langle a,b\rangle\in A^2: a\leftrightarrow b
\in F\}
;\\
& \theta\,\mapsto\,\{a\in A:
a\wedge e\equiv_\theta e
\}.
\end{eqnarray*}
We abbreviate $\sbA/\leibniz^\sbA F$ as $\sbA/F$, and $a/\leibniz^\sbA F$ as $a/F$, noting that
\[
\textup{$a\rig b\in F$ iff
$a/F\leqslant b/F$
in $\sbA/F$.}
\]
Whenever $\sbB$ is a subalgebra of $\sbA$, and $F$ is a deductive filter of $\sbA$, then $B\cap F$
is a deductive filter of $\sbB$ and
\begin{equation}\label{omega and subalgebras}
\leibniz^\sbB(B\cap F)=
(\leibniz^\sbA F)|_B.
\end{equation}

Because of the above isomorphism,
$\sbA$ is FSI
iff its smallest deductive filter $\text{$\{a\in A\colon e\leqslant a\}$}$ is meet-irreducible in its lattice of deductive
filters, and that amounts to the join-irreducibility of $e$ in $\langle A;\wedge,\vee\rangle$.
Since this last condition is expressible as a universal first order sentence,
every variety of S[I]RLs has EDPM.
Therefore, Theorem~\ref{thm:ESWitnessedByFSI} applies to all such
varieties.

For any subset $X$ of
$\sbA$, the smallest deductive filter of $\sbA$ containing $X$ is
denoted by $\Fg^{\alg{A}}X$.  Thus, $\Fg^\sbA X$ consists of all $a\in A$ such that
\[
a\geqslant x_1\wedge\,\dots\,\wedge x_n \textup{ for some }x_1,\dots,x_n\in X\cup\{e\}, \textup{ where } 0<n\in\omega.
\]
In particular, if $e\geqslant b\in A$, then $\Fg^\sbA\{b\}=\{a\in A:b\leqslant a\}$.

It follows that the deductive filters of a subalgebra $\alg{B}$ of $\alg{A}$ are just the sets $B\cap F$ such that $F$ is a deductive filter of $\sbA$.
Note that the deductive filters of a Brouwerian algebra are exactly its
(lattice-) filters.

A filter $F$ of a lattice $\langle L;\wedge,\vee\rangle$
is said to be \emph{prime} if its complement $L\ld F$ is closed under the binary operation $\vee$.
For any filters $F,G,H$ of $\langle L;\wedge,\vee\rangle$,
\begin{equation}\label{lem:filterProperties:ForG}
\textup{if $F \cap G \subseteq H$ and $H$ is prime, then $F \subseteq H$ or $G \subseteq H$.}
\end{equation}
The \emph{Prime Filter
Extension Theorem} asserts that, when
$\langle K;\wedge,\vee\rangle$ is a sublattice of a distributive lattice $\langle L;\wedge,\vee\rangle$, then the prime filters of
$\langle K;\wedge,\vee\rangle$ are exactly
the non-empty sets $K\cap F$ such that $F$ is a prime filter
of $\langle L;\wedge,\vee\rangle$ \cite[Theorem~III.6.5]{BD74}.

We use
$\Pr(\alg{A})$ to denote
the set of all prime deductive filters of an S[I]RL
$\alg{A}$, \emph{including} $A$ itself.  We always consider $\Pr(\alg{A})$ to be partially ordered by set inclusion.
For a
deductive
filter $F$
of $\alg{A}$,
we write
$$\upclose^{\alg{A}} F = \{ H \in \Pr(\alg{A}) : F \subseteq H \},$$
i.e., $\upclose^{\alg{A}} F$ abbreviates $\upclose^{\,\Pr(\alg{A})}F$, where $\Pr(\alg{A})$ is considered as a subset
of the lattice of deductive filters of $\sbA$.
\begin{remark}\label{correspondence theorem}
Suppose $h\colon\sbA\mrig\sbB$ is a surjective homomorphism of S[I]RLs.  The kernel of $h$ is
$\leibniz^\sbA K$ for some deductive filter $K$ of $\sbA$.
If $G$ is a deductive filter of $\sbA$, with $K\subseteq G$, then $h[G]\seteq\{h(g):g\in G\}$
is a deductive filter of $\sbB$, and by
the Correspondence Theorem of Universal Algebra,
\[
H\mapsto
h^{\leftarrow}[H]\seteq\{a\in A:h(a)\in H\}
\]
is a lattice isomorphism from the deductive filter lattice of
$\alg{B}$ onto the lattice of deductive filters of $\sbA$ that contain $K$; the inverse
isomorphism is given by $G\mapsto h[G]$.  In particular,
\begin{equation}\label{lem:filterProperties:class}
h^\leftarrow [h[G]]=G \textup{ for all deductive filters $G$ of $\sbA$ such that $K\subseteq G$.}
\end{equation}
Clearly, a deductive filter
$H$ of $\alg{B}$ is prime iff the filter $h^\leftarrow[H]$ of $\sbA$ is prime, so $H\mapsto
h^{\leftarrow}[H]$ also defines an isomorphism of partially ordered sets from $\Pr(\alg{B})$ onto $\upclose^{\alg{A}}K$
(both ordered by inclusion).
\end{remark}

\section{Negative Cones}

\begin{definition}\label{defn:negative cone}
The \emph{negative cone}
of an S[I]RL $\alg{A} = \langle A; \wedge, \vee,\bdot, \to, e \,[, \neg] \rangle$
is the Brouwerian algebra
$$\alg{A}^-=\left\langle A^-;\, \wedge|_{(A^-)^2},\, \vee|_{(A^-)^2},\, \to^-,\, e \right\rangle\!\textup{,}$$
where $A^- = \{ a \in A : a \leqslant e \}$ and $a \to^- b = (a \to b) \wedge e$ for all $a,b \in A^-$.
\end{definition}

\begin{lemma}
\label{lem:dualityProperties}
Let $\alg{A}$ and\/ $\alg{B}$ be S[I]RLs, and $F$ a deductive filter of\/ $\alg{A}$\textup{.}
\begin{enumerate}
\item \label{lem:dualityProperties:minus} If\/ $h \colon \alg{A} \mrig \alg{B}$ is a
homomorphism, then $h|_{A^-}$ is a
homomorphism from $\alg{A}^-$ to $\alg{B}^-$\textup{.} If, moreover, $h$ is surjective, then so is\/ $h|_{A^-}$\textup{.}
\item \label{lem:dualityProperties:FilterRestriction}
$A^-\cap F$ is a filter of $\alg{A}^-$\textup{,} and\/ $(\alg{A}/F)^- \cong \alg{A}^-/(A^-\cap F)$\textup{,} the
isomorphism and its inverse
being
\[
\textup{$a/F \mapsto (a \wedge e)/(A^-\cap F)$ and\/ $a/(A^-\cap F)\mapsto a/F$}\textup{.}
\]
\end{enumerate}
\end{lemma}

\begin{proof}
(\ref{lem:dualityProperties:minus})
Since homomorphisms between S[I]RLs are isotone maps that preserve $e$\textup{,} we have
$h[A^-] \subseteq B^-$.  Also, as $h$ is a homomorphism, it is clear from the definitions of the operations on the negative cone
that $h|_{A^-}$ is a homomorphism from $\sbA^-$ to $\sbB^-$.
Now suppose $h$ is onto. For each $b \in B^-$, there exists $a \in A$ with $h(a) = b$, and since $b \leqslant e$, we have
$ b = h(a) \wedge e = h(a \wedge e)$.  As
$a \wedge e \in A^-$, this shows that $B^-=h[A^-]$.

(\ref{lem:dualityProperties:FilterRestriction})
Clearly, $A^-\cap F$ is a filter of $\sbA^-$.
Let $q \colon \alg{A} \mrig \alg{A}/F$ be the canonical surjection. By (\ref{lem:dualityProperties:minus}), $q|_{A^-} \colon \alg{A}^- \mrig (\alg{A}/F)^-$ is a surjective
homomorphism.  For all $a,b\in A^-$, we have $a\leftrightarrow b\in F$ iff $(a\leftrightarrow b)\wedge e\in A^-\cap F$ (since
$F$ is upward closed and contains $e$), i.e.,
the kernel of $q|_{A^-}$ is $\leibniz^{\alg{A}^-}\!(A^-\cap F)$.  Thus, by the Homomorphism Theorem,
$a/(A^-\cap F)\mapsto a/F$ defines an isomorphism from $\sbA^-/(A^-\cap F)$ onto $(\alg{A}/F)^-$.
For any $a\in A$, if $a/F \in (A/F)^-$, then $a\wedge e\in A^-$ and
$(a \wedge e)/ F = (a/F) \wedge (e/F) = a/F$,
so
$a/F \mapsto (a \wedge e)/(A^-\cap F)$ defines the inverse of the above isomorphism.
\end{proof}

Given an S[I]RL $\sbA$, if $F$ is a filter of $\sbA^-$, then
\begin{equation}\label{going up}
\textup{$\Fg^\sbA\!F=\{a\in A:a\geqslant b\textup{ for some }b\in F\}$, \,so\, $A^-\cap\Fg^\sbA\!F=F$.}
\end{equation}

\begin{definition}
An S[I]RL $\sbA$ is
\emph{negatively generated} if
$A=\Sg^\sbA (A^-)$.
\end{definition}
As surjective homomorphisms always map generating sets onto generating sets, the following
lemma applies.

\begin{lemma}\label{neg cone gen in images}
If\/ $h\colon\sbA\mrig\sbB$ is a surjective homomorphism of S[I]RLs and\/ $\sbA$ is negatively generated,
then so is\/ $\sbB$\textup{.}
\end{lemma}

\section{Duality for Brouwerian Algebras}

A well-known duality between Heyting algebras and Esakia spaces was established in \cite{Esa74,Esa85}.
It entails an analogous duality between the variety $\mathsf{BRA}$ of Brouwerian algebras
(considered as a concrete category, equipped with all algebraic homomorphisms) and the category
$\mathsf{PESP}$ of `pointed Esakia spaces' defined below, i.e., there is a category equivalence between
$\mathsf{BRA}$ and the opposite category of $\mathsf{PESP}$.  This is explained, for instance, in \cite[Section~3]{BMR17}, but
we recall the key definitions here.

A structure $\spc{X}=\langle X; \tau, \leqslant, m \rangle$ is a \emph{pointed Esakia space} if
$\langle X; \leqslant \rangle$ is a partially ordered set with a greatest element $m$, and $\langle X;\tau \rangle$
is a compact Hausdorff space in which
\begin{enumerate}
\renewcommand{\theenumi}{\roman{enumi}}
\item every open set is a union of clopen (i.e., closed and open) sets;
\item $\upclose x$
is closed, for all $x \in X$;
\item $\downclose V$
is clopen, for all clopen $V \subseteq X$.
\end{enumerate}
In this case, the {\em Priestley separation axiom\/} of \cite{Pri70} holds: for all $x,y \in X$,
\begin{enumerate}
\renewcommand{\theenumi}{\roman{enumi}}
\setcounter{enumi}{3}
\item if $x \nleqslant y$, then $x \in U$ and $y \notin U$, for some clopen up-set $U \subseteq X$.
\end{enumerate}
The morphisms of $\cl{PESP}$ are the so-called \emph{Esakia morphisms} between
these objects.  They are the isotone continuous functions $g \colon \spc{X} \mrig \spc{Y}$
such that,
\begin{equation}\label{esakia morphism}
\textup{whenever $x \in X$ and $g(x) \leqslant y \in Y$, then $y = g(z)$ for some $z \in \upclose x$.}
\end{equation}
It follows that $g(m) = m$, and if $g$ is bijective then $g^{-1} \colon\spc{Y} \mrig \spc{X}$ is also an Esakia morphism, so $g$ is a (categorical) isomorphism.

Given $\alg{A}\in\mathsf{BRA}$ and $a\in A$,
let $\varphi^{\alg{A}}(a)$ denote
$\{F\in\Pr(\alg{A}): a\in F\}$ and $\varphi^{\alg{A}}(a)^c$ its complement $\{F\in\Pr(\alg{A}): a\notin F\}$.
The \emph{dual} (in $\cl{PESP}$) of
$\alg{A}$ is $\alg{A}_* = \langle \Pr(\alg{A}); \tau, \subseteq, A \rangle$, where $\tau$ is the topology on $\Pr(\sbA)$ with sub-basis
\[
\{ \varphi^{\alg{A}}(a) : a \in A \} \cup \{ \varphi^{\alg{A}}(a)^c : a \in A \}.
\]

The \emph{dual} of a morphism $h \colon \alg{A} \mrig \alg{B}$ in $\mathsf{BRA}$ is the Esakia morphism
$h_* \colon \alg{B}_* \mrig \alg{A}_*$, defined by $F \mapsto h^{\leftarrow}[F]$.
Thus, the contravariant functor $(-)_* \colon \cl{BRA} \mrig \cl{PESP}$ is given by $\alg{A} \mapsto \alg{A}_*$; $h \mapsto h_*$.

The contravariant functor $(-)^* \colon \cl{PESP} \mrig \cl{BRA}$
works as follows.  Let $g \colon \spc{X} \mrig \spc{Y}$ be an Esakia morphism, where $\spc{X}, \spc{Y} \in \cl{PESP}$. Then
\[
\spc{X}^* = \langle \Cpu(\spc{X}); \cap, \cup, \to, X \rangle \in \cl{BRA},
\]
where
$\Cpu(\alg{X})$ is the set of all \emph{non-empty} clopen up-sets of $\spc{X}$, and
\[
U \to V \seteq X \setminus \downclose(U \setminus V)
\]
for all $U, V \in \Cpu(\alg{X})$.
The homomorphism $g^* \colon \spc{Y}^* \mrig \spc{X}^*$ is given by $U \mapsto g^\leftarrow[U]$.  We refer to $\spc{X}^*$ [resp.\ $g^*$] as the \emph{dual}
of $\sbX$ [resp.\ $g$] in $\mathsf{BRA}$.

The functors $(-)_*$ and $(-)^*$ establish the aforementioned duality.
For $\alg{A}\in\class{BRA}$ and $\spc{X}\in\class{PESP}$, the respective canonical isomorphisms from $\alg{A}$ to ${\alg{A}_*}^*$
and from $\spc{X}$ to ${\spc{X}^*}_*$ are given by $a\mapsto\varphi^{\alg{A}}(a)$ and
$$x \mapsto \{ U \in \Cpu(\spc{X}) : x \in U \}.$$

In $\cl{PESP}$, there is a notion of substructure: an \emph{E-subspace} of $\spc{X} \in \cl{PESP}$ is a non-empty closed up-set $U$ of $\spc{X}$.  It is the universe of a pointed Esakia space $\spc{U}$, with the restricted order and the subspace topology, so the inclusion $\spc{U} \mrig \spc{X}$ is an Esakia morphism.

\begin{lemma}
\label{lem:BRAProperties}
\textup{(\cite{Esa85})}
\begin{enumerate}
\item \label{lem:BRAProperties:SurjectiveVSInjective} A homomorphism $h$ between Brouwerian algebras is surjective iff\/ $h_*$ is injective. Also, $h$ is injective iff\/ $h_*$ is surjective.
\item \label{lem:BRAProperties:Image} The image of a morphism in $\cl{PESP}$ is an $E$-subspace of the co-domain.
\end{enumerate}
\end{lemma}

The ES property for $\mathsf{BRA}$ is relied on in the right-to-left implication of the first assertion of Lemma~\ref{lem:BRAProperties}(\ref{lem:BRAProperties:SurjectiveVSInjective}).
The forward implication in the second assertion of (\ref{lem:BRAProperties:SurjectiveVSInjective}) employs the Prime Filter
Extension Theorem.  In the absence of a convenient reference, a proof of the next lemma is supplied below; the result is
presumably well-known.

\begin{lemma}
\label{lem:subspace}
Let\/ $F$ be a filter of a Brouwerian algebra\/ $\alg{A}
$\textup{,} and\/
${q\colon \alg{A} \mrig \alg{A}/F}$
the canonical surjection. Then $q_*$ is an isomorphism from $(\alg{A}/F)_*$ onto the E-subspace $q_*[(\alg{A}/F)_*]$ of $\alg{A}_*$
whose universe is\/ $\upclose^\alg{A} F$.
Also, the map
$$\varphi^{\alg{A}}_F \colon a/F \mapsto \{ H \in \Pr(\alg{A}) : F\cup\{a\}\subseteq H \}$$
is an isomorphism from $\alg{A}/F$ onto $(q_*[(\alg{A}/F)_*])^*$ and the following diagram commutes, where $i_1 \colon q_*[(\alg{A}/F)_*] \mrig \alg{A}_*$ is the inclusion map.
\begin{center}
\begin{picture}(125,65)

\put(5,5){${\alg{A}_*}^*$}

\put(25,5){$\xrightarrow{\hspace{30pt}}$}
\put(25,5){$\xrightarrow{\hspace{28pt}}$}
\put(37,13){\tiny $i_1^*$}

\put(65,5){$(q_*[(\alg{A}/F)_*])^*=(\upclose^{\alg{A}}F)^*$}

\put(8,45){$\alg{A}$}

\put(20,45){$\xrightarrow{\hspace{45pt}}$}
\put(20,45){$\xrightarrow{\hspace{43pt}}$}
\put(43,53){\tiny $q$}

\put(75,45){$\alg{A}/F$}

\put(10,40){\rotatebox{270}{$\xrightarrow{\hspace{15pt}}$}}
\put(7,37){\rotatebox{270}{\scalebox{1.5}[1]{$\sim$}}}
\put(16,29){\tiny $\varphi^\alg{A}$}

\put(87,40){\rotatebox{270}{$\xrightarrow{\hspace{15pt}}$}}
\put(84,37){\rotatebox{270}{\scalebox{1.5}[1]{$\sim$}}}
\put(93,29){\tiny $\varphi^{\alg{A}}_F$}

\end{picture}
\end{center}
Furthermore, if $G$ is a
filter of $\alg{A}$\textup{,} with
$F \subseteq G$\textup{,} then the following diagram commutes, where $q' \colon a/F \mapsto a/G$\textup{,} and $i_2$ is the inclusion map.
\begin{center}
\begin{picture}(160,65)


\put(-52,5){$(\upclose^{\alg{A}}F)^*=(q_*[(\alg{A}/F)_*])^*$}

\put(65,5){$\xrightarrow{\hspace{20pt}}$}
\put(65,5){$\xrightarrow{\hspace{18pt}}$}
\put(75,13){\tiny $i_2^*$}

\put(95,5){$(q_*[(\alg{A}/G)_*])^*=(\upclose^{\alg{A}}G)^*$}

\put(25,45){$\alg{A}/F$}

\put(52,45){$\xrightarrow{\hspace{45pt}}$}
\put(52,45){$\xrightarrow{\hspace{43pt}}$}
\put(75,53){\tiny $q'$}

\put(105,45){$\alg{A}/G$}

\put(30,40){\rotatebox{270}{$\xrightarrow{\hspace{15pt}}$}}
\put(27,37){\rotatebox{270}{\scalebox{1.5}[1]{$\sim$}}}
\put(36,29){\tiny $\varphi^{\alg{A}}_F$}

\put(117,40){\rotatebox{270}{$\xrightarrow{\hspace{15pt}}$}}
\put(114,37){\rotatebox{270}{\scalebox{1.5}[1]{$\sim$}}}
\put(123,29){\tiny $\varphi^{\alg{A}}_G$}

\end{picture}
\end{center}
\end{lemma}

\begin{proof}
As $q \colon \alg{A} \mrig \alg{A}/F$ is a surjective $\mathsf{BRA}$-morphism,
Lemma~\ref{lem:BRAProperties}(\ref{lem:BRAProperties:SurjectiveVSInjective})
shows that
$q_*\colon (\alg{A}/F)_* \mrig \alg{A}_*$ is an injective Esakia morphism; its image
is the E-subspace $\upclose^\sbA F$ of $\sbA_*$, by Lemma~\ref{lem:BRAProperties}(\ref{lem:BRAProperties:Image}) and
Remark~\ref{correspondence theorem}.
Let $k= {q_*}^{-1}\big|_{\upclose^\sbA F}$\,,
so $k\colon
\,\upclose^\sbA F\cong (\sbA/F)_*$
is defined by
\[
k(H)=q[H]= H/F\seteq\{a/F:a\in H\} \ \ (H\in\upclose^\sbA F),
\]
and $k^*\colon{(\sbA/F)_*}^*\cong (\upclose^\sbA F)^*$.  Since $\varphi^{\sbA/F}\colon\sbA/F\cong {(\sbA/F)_*}^*$, we have
\[
k^*\circ\varphi^{\sbA/F}\colon\sbA/F\cong(\upclose^\sbA F)^*.
\]
For each $a\in A$,
\begin{align*}
& (k^*\circ\varphi^{\sbA/F})(a/F)\,=\,k^*(\{H\in \Pr(\sbA/F):a/F\in H\})\\
& \quad\quad\quad\quad\quad\quad\quad\;\;\: =\,\{G\in \upclose^\sbA F:a/F\in G/F\}\\
& \quad\quad\quad\quad\quad\quad\quad\;\;\: =\,\{G\in \upclose^\sbA F:a\in G\} \textup{ \,(by (\ref{lem:filterProperties:class}))} = \varphi^\sbA_F(a/F),
\end{align*}
so
$k^*\circ\varphi^{\sbA/F}=\varphi^\sbA_F$, whence $\varphi^\sbA_F\colon\sbA/F\cong(\upclose^\sbA F)^*=
(q_*[(\alg{A}/F)_*])^*$.

Commutativity of the second diagram subsumes that of the first (after identification of $\sbA/\{e\}$ with $\sbA$, and
$\varphi^\sbA_{\{e\}}$ with $\varphi^\sbA$).

Accordingly, let $G$ be a filter of $\sbA$, with $F \subseteq G$, so $q' \colon a/F \mapsto a/G$ is a homomorphism from
$\sbA/F$ onto $\sbA/G$.  For each $a\in A$, the respective left and right hand sides of the equation
$\varphi^\sbA_G(q'(a/F))=i_2^*(\varphi^\sbA_F(a/F))$
are, by definition,
\[
\textup{$\{H\in\Pr(\alg{A}):G\cup\{a\}\subseteq H\}$ \,and\,
$(\upclose^\sbA G)\cap\{H\in\Pr(\alg{A}):F\cup\{a\}\subseteq H\}$,}
\]
which are clearly equal, so $i_2^* \circ \varphi^{\alg{A}}_F = \varphi^{\alg{A}}_G \circ q'$.
\end{proof}

\section{Depth}

For $\spc{X} = \langle X; \tau, \leqslant, m \rangle \in \cl{PESP}$ and $x \in X$, we define $\depth{(x)}$ (the \emph{depth} of $x$ in $\spc{X}$) to be the greatest $n \in \omega$ (if it exists) such that there is a chain
$$ x = x_0 < x_1 < \dots < x_n = m$$
in $\spc{X}$.  (Thus, $m$ has depth $0$ in $\spc{X}$.) If no greatest such $n$ exists, we set $\depth{(x)}=\infty$.
We define $\depth(\spc{X}) = \sup\,\{\depth(x) : x \in X\}$.

For $\alg{A} \in \cl{BRA}$ and any subvariety $\cl{K}$ of $\cl{BRA}$, we define
\[
\textup{$\depth(\alg{A}) = \depth(\alg{A}_*)$ \,and\, $\depth(\cl{K}) = \sup\,\{\depth(\alg{B}) : \alg{B} \in \cl{K}\}$.}
\]
If a subvariety of $\cl{BRA}$ is finitely generated (i.e., of the form $\mathbb{V}(\sbA)$ for some finite $\sbA\in\cl{BRA}$)
then it has finite depth, and if it has finite depth then it is
locally finite (i.e., its finitely generated members are finite).  Both converses are false.
For each $n \in \omega$, the class
${\cl{BRA}_n = \{ \alg{A} \in \cl{BRA} : \depth(\alg{A}) \leq n \}}$
is a
finitely axiomatizable variety.
These claims are explained in \cite[Section~4]{BMR17},
where their antecedents are also discussed.

\begin{definition}
\label{def:depth}
For any S[I]RL $\sbA$ and any variety $\cl{K}$ of S[I]RLs, we define the
\emph{depth} of
$\alg{A}$ to be
the depth of its negative cone $\alg{A}^-$, and the \emph{depth} of $\cl{K}$ to be
$\sup\,\{\depth(\alg{B}) : \alg{B} \in \cl{K}\}$.
\end{definition}

For each $n\in\omega$, an S[I]RL
has depth at most $n$
iff it satisfies the equations that result from the axioms for $\cl{BRA}_n$
when we replace $\rig$ by $\rig^-$ and $x$ by $x\wedge e$, for every apparent variable $x$.  Thus, the S[I]RLs of depth at most $n$
also form a finitely axiomatizable variety.

Consequently, every finitely generated variety of S[I]RLs has
finite depth, and the class of S[I]RLs of depth greater than (any fixed) $n$ is closed
under ultraproducts.  Since the class of S[I]RLs that are FSI is also closed under ultraproducts, it follows that,
when each FSI member of a variety $\cl{K}$ of S[I]RLs has finite depth, then so does $\cl{K}$.

\section{The ES Property}

We can now formulate
the main result of this paper.

\begin{theorem}
\label{thm:SIRLhasES}
Let\/ $\cl{K}$ be a variety of S[I]RLs,
such that each FSI member of\/ $\cl{K}$
has finite depth and is negatively generated.
Then every\/ $\cl{K}$-epimorphism is surjective.
\end{theorem}

The proof of Theorem~\ref{thm:SIRLhasES} is by contradiction, and it proceeds via a sequence of claims.
Let $\cl{K}$ be as postulated, and suppose
that $\cl{K}$ lacks the ES property.
By Theorem~\ref{thm:ESWitnessedByFSI}, some $\alg{A} \in \cl{K}_{FSI}$ has a proper $\cl{K}$-epic subalgebra.
Now $\sbA/\theta\in\cl{K}$ for all $\theta\in\mathit{Con}\,\alg{A}$, as $\cl{K}$ is a variety.
We shall define a congruence $\theta$ of $\alg{A}$ such that
the following is true.

\smallskip

\begin{claim}
\label{claim:subalgebra}
\emph{There exist\/ $a \in A$ and a\/ $\cl{K}$-epic proper subalgebra\/ $\alg{C}$ of\/ $\alg{A}/\theta$ such that\/
$\alg{C}$ is negatively generated
and\/
$\sbA/\theta$ is generated by\/ $C^-\cup\{a/\theta\}$\textup{,} and\/
$a/\theta \prec e/\theta$ in\/ $\alg{A}/\theta$\textup{.}}
\end{claim}

\smallskip

Once $\theta$ has been identified and Claim~\ref{claim:subalgebra} proved, we shall
contradict the fact that $\sbC$ is $\cl{K}$-epic in
$\sbA/\theta$, by constructing a
non-identity homomorphism $\ell \colon\alg{A}/\theta \mrig \alg{A}/\theta$, such that $\ell|_C = \id_C$,
as follows.

Let $b=a/\theta$.  Then $\Fg^{\sbA/\theta}\{b\}=\{d\in A/\theta:b\leqslant d\}$, because $b<e/\theta$.
Let $\alpha = \leibniz^{\alg{A}/\theta}{\,\Fg^{\sbA/\theta}\{b\}}$.  For any $u,v\in A/\theta$, we have
\begin{equation}\label{uv}
\textup{$u\equiv_\alpha v$ \,iff\, $b\leqslant u\leftrightarrow v$.}
\end{equation}
In particular, $e/\theta\equiv_\alpha b$, by (\ref{t order}) and (\ref{t laws}).

Let $\{a_i : i \in I\}$ be an indexing of $A/\theta$, and
$\vec{c} = c_0, c_1, \dots$ a well-ordering of the elements of $C^-$.
Since $C\neq A/\theta$ and $\sbA/\theta$ is generated by $C^-\cup\{b\}$, it follows that
$b \notin C$ and,
for each $i \in I$, we have $a_i=t_i^{\alg{A}/\theta}(b,\vec{c})$ for a suitable
S[I]RL-term $t_i(x,\vec{y})$, where $\vec{y} = y_0, y_1, \dots\,\,$.
When $a_i$ is some $c_j \in C^-$, we can (and do) choose $t_i$ to be the variable $y_j$.

We define $\ell\colon A/\theta \mrig A/\theta$ by
\[
\textup{$\ell(a_i) = t_i^{\alg{A}/\theta}(e/\theta,\vec{c})$ \,for all $i\in I$.}
\]
By the above choice, $\ell(c_j) = c_j$ for $j = 0, 1, \dots$,
while $\ell[A/\theta] \subseteq C$, because $\sbC$ is a subalgebra of $\sbA/\theta$.
We claim that $\ell$ is a homomorphism.

To see this,
let $\sigma$ be a
fundamental S[I]RL-operation symbol, and let
$a_{i_1}, \dots, a_{i_n} \in A/\theta$, where $n$ is the rank of $\sigma$.
Then $\sigma^{\alg{A}/\theta}(a_{i_1}, \dots, a_{i_n}) = a_j$ for some $j \in I$.
For this $j$, we perform the following
calculation, where every term is evaluated in $\alg{A}/\theta$\,:
\begin{align*}
\sigma(\ell(a_{i_1}), \dots, \ell(a_{i_n})) &= \sigma(t_{i_1}(e/\theta,\vec{c}), \dots, t_{i_n}(e/\theta,\vec{c})) \\
&\equiv_\alpha \sigma(t_{i_1}(b,\vec{c}), \dots, t_{i_n}(b,\vec{c})) \\
&= \sigma(a_{i_1}, \dots, a_{i_n})
=a_j = t_j(b,\vec{c}) \\
&\equiv_\alpha t_j(e/\theta,\vec{c}) = \ell(a_j) = \ell(\sigma(a_{i_1}, \dots,a_{i_n})).
\end{align*}
By (\ref{uv}), therefore,
\[
b \leqslant (\sigma(\ell(a_{i_1}), \dots, \ell(a_{i_n})) \leftrightarrow \ell(\sigma(a_{i_1}, \dots,a_{i_n})))\wedge (e/\theta).
\]
Note that
\[
(\sigma(\ell(a_{i_1}), \dots, \ell(a_{i_n})) \leftrightarrow \ell(\sigma(a_{i_1}, \dots,a_{i_n})))\wedge (e/\theta) \in C
\]
(because ${\ell[A/\theta] \subseteq C}$), but $b\notin C$, so
\[
b < (\sigma(\ell(a_{i_1}), \dots, \ell(a_{i_n})) \leftrightarrow \ell(\sigma(a_{i_1}, \dots,a_{i_n})))\wedge (e/\theta)\leqslant e/\theta.
\]
Since $b\prec e/\theta$ in $\sbA/\theta$, this forces
\[
e/\theta=(\sigma(\ell(a_{i_1}), \dots, \ell(a_{i_n})) \leftrightarrow \ell(\sigma(a_{i_1}, \dots,a_{i_n})))\wedge (e/\theta),
\]
i.e., $e/\theta\leqslant \sigma(\ell(a_{i_1}), \dots, \ell(a_{i_n})) \leftrightarrow \ell(\sigma(a_{i_1}, \dots,a_{i_n}))$.
Then, by (\ref{t reg}),
\[
\sigma(\ell(a_{i_1}), \dots, \ell(a_{i_n})) = \ell(\sigma(a_{i_1}, \dots,a_{i_n})),
\]
confirming that $\ell$ is a homomorphism.

For each $c \in C$, we have $c = t^{\alg{A}/\theta}(\vec{c})$ for some S[I]RL-term $t$ (as $\alg{C}$ is generated by $C^-$), so
\[
\ell(c)=\ell(t(\vec{c}))= t(\ell(c_0),\ell(c_1), \dots)= t(c_0,c_1, \dots)=c.
\]
This shows that
$\ell|_C = \id_C$,
but $\ell(b) \neq b$, since $\ell(b) \in C$.  As intended, this contradicts the fact that
$\sbC$ is $\cl{K}$-epic in $\sbA/\theta$.

\medskip

It remains to construct $\theta$ and to prove Claim~\ref{claim:subalgebra}.  The construction of $\theta$
exploits the assumption that members of $\cl{K}_{FSI}$ have finite depth.

Recall that $\alg{A}$
has a proper $\cl{K}$-epic subalgebra,
$\alg{B}$ say.  As $\sbA$ is FSI, so is $\alg{B}$ (because $\cl{K}$ has EDPM, as we noted after (\ref{omega and subalgebras})).
By assumption, therefore, $\sbA$ and $\sbB$ are both negatively generated,
so
$\alg{B}^- \neq \alg{A}^-$, because $\sbB\neq \sbA$.

For each $F\in\Pr(\sbA^-)$, we clearly have $B\cap F=B^-\cap F=i_*(F)$, where
$i$ is the inclusion map $i \colon \alg{B}^- \mrig \alg{A}^-$, considered as a $\cl{BRA}$-morphism.

As $i$ is
not surjective, its dual
$i_* \colon (\alg{A}^-)_* \mrig (\alg{B}^-)_*$ is not injective, by Lemma~\ref{lem:BRAProperties}(\ref{lem:BRAProperties:SurjectiveVSInjective}), i.e., the following set is not empty:
\[
W\seteq \{\langle F_1,F_2\rangle\in (\Pr(\alg{A}^-))^2: F_1 \neq F_2 \text{ and } F_1 \cap B = F_2 \cap B\}.
\]
By assumption, $\sbA^-$ has finite depth, so
\begin{equation*}
\left\{
\min\,\{
\depth(F_1),\depth(F_2)\}: \langle F_1,F_2\rangle\in W\right\}
\end{equation*}
is a non-empty subset of $\omega$,
and therefore has a least element, $n$ say.  Pick $F_1\in\Pr(\alg{A}^-)$ such that $\depth(F_1) = n$
and $\langle F_1,G\rangle\in W$ for some $G$.  Now,
\begin{equation}\label{eq:DefOfF1}
\textup{whenever $\langle F_1',F_2'\rangle\in W$, then $\depth(F_1) \leq \depth(F_1'), \depth(F_2')$.}
\end{equation}

Having fixed $F_1$ in this way, we similarly choose $F_2 \in \Pr(\alg{A}^-) \setminus \{F_1\}$ such that
$F_1 \cap B = F_2 \cap B$ and
\begin{equation}\label{eq:DefOfF2}
\textup{whenever $\langle F_1,F_2'\rangle\in W$, then $\depth(F_2) \leq \depth(F_2')$.}
\end{equation}
As $\langle F_1,F_2\rangle\in W$, we have $\depth(F_1) \leq \depth(F_2)$ (by (\ref{eq:DefOfF1})), so $F_1$ is not a proper
subset of $F_2$.

\begin{lemma}
\label{lem:F1ltG}
If\/ $F_1\subsetneq G \in \Pr(\alg{A}^-)$\textup{,}
then\/ $F_2 \subsetneq  G$\textup{.}
\end{lemma}

\begin{proof}
Let $F_1\subsetneq G \in \Pr(\alg{A}^-)$, so $\depth(G) < \depth(F_1)$.  As
$i_*$ is an Esakia morphism and
$i_*(F_2) = i_*(F_1) \subseteq
i_*(G)$, there exists $H \in \Pr(\alg{A}^-)$ such that $F_2 \subseteq
H$ and $i_*(G) = i_*(H)$, by (\ref{esakia morphism}),
i.e., $G\cap B= H\cap B$.  Therefore, if ${G\neq H}$, then
$\depth(F_1) \leq \depth(G)$, by (\ref{eq:DefOfF1}).  This is a contradiction, so
$G = H$, whence $F_2 \subseteq
G$.  If $F_2 = G$, then
\[
\textup{$\depth(F_2) = \depth(G) < \depth(F_1)$,}
\]
contradicting the fact that
$\depth(F_1) \leq \depth(F_2)$. Therefore, $F_2 \subsetneq  G$.
\end{proof}

\begin{lemma}
\label{lem:F2ltG}
If\/ $F_2\subsetneq G \in \Pr(\alg{A}^-)$\textup{,}
then\/ $F_1 \subseteq
G$\textup{.}
\end{lemma}

\begin{proof}
Let $F_2\subsetneq G \in \Pr(\alg{A}^-)$.
Again, $i_*(F_1) = i_*(F_2) \subseteq
i_*(G)$, so there exists $H \in \Pr(\alg{A}^-)$ such that $F_1 \subseteq
H$ and $i_*(G)=i_*(H)$. Suppose $G \neq H$. If $F_1 = H$, then $F_1 \cap B = G \cap B$, so, by (\ref{eq:DefOfF2}), $\depth(F_2) \leq \depth(G)$, contradicting the fact that $F_2\subsetneq G$.
Therefore, $F_1 \subsetneq  H$, so $\depth(H) < \depth(F_1)$. Then, by (\ref{eq:DefOfF1}), $\depth(F_1) \leq \depth(H)$, since $G \cap B = H \cap B$.  This is a contradiction, so $G = H$, whence $F_1 \subseteq G$.
\end{proof}

Recalling that $F_1,F_2$ are distinct and that $F_1$ is not properly contained in $F_2$,
we make the following claim:

\smallskip

\begin{claim}
\label{claim:2possibilities}
\emph{There are just two possibilities:}
\begin{enumerate}
\renewcommand{\theenumi}{\Alph{enumi}}
\item \label{claim:2possibilities:F2ltF1}
\emph{$F_2 \subsetneq  F_1$\textup{,} in which case $F_2 \prec F_1$ (in fact, $F_1$ is the least strict upper bound
of $F_2$ in\/ $\Pr(\sbA^-)$\textup{);}}
\item \label{claim:2possibilities:incomparable}
\emph{$F_1$ and $F_2$ are incomparable, in which case they have the same depth,
the same strict upper bounds and, therefore, the same covers in\/ $\Pr(\sbA^-)$.}
\end{enumerate}
\end{claim}

\begin{proof}
If $F_2\subsetneq F_1$, then $F_1$ is the least strict upper bound of $F_2$,
by Lemma~\ref{lem:F2ltG}, i.e., $F_1$ is the unique cover of $F_2$.
We may therefore assume that $F_2$ is not a proper subset of $
F_1$, i.e., that $F_1$ and $F_2$ are incomparable.  Then,
by Lemmas~\ref{lem:F1ltG} and \ref{lem:F2ltG},
$F_1$ and $F_2$ have the same strict upper bounds, and hence the same covers in $\Pr(\sbA^-)$, from
which it clearly follows that they also have the same depth.
\end{proof}

By (\ref{going up}), $\Fg^{\alg{A}}(F_1 \cap F_2)=\{d\in A:d\geqslant c \textup{ for some }c\in F_1\cap F_2\}$.
We define
\[
\theta = \leibniz^{\alg{A}}\,{\Fg^{\alg{A}}(F_1 \cap F_2)}.
\]

\begin{claim}
\label{claim:retract}
\emph{The following diagram commutes, where the maps will be defined in the proof.}
\begin{center}
\begin{picture}(200,105)


\put(25,5){$\spc{Z}^*$}

\put(40,5){$\xrightarrow{\hspace{45pt}}$}
\put(48,1){\scalebox{4}[1]{$\sim$}}
\put(57,13){\tiny $(i_*|_Y)^*$}

\put(95,5){$\spc{Y}^*$}

\put(113,12){\rotatebox{180}{$\xrightarrow{\hspace{25pt}}$}}
\put(119,3){\scalebox{2.5}[1]{$\sim$}}
\put(124,15){\tiny $\varphi^{\alg{A}^-}_{F_1}$}

\put(150,5){$\alg{A}^-\!/F_1$}

\put(5,45){$\alg{B}^-\!/(F_1 \cap B)$}

\put(95,45){$\spc{X}^*$}

\put(113,52){\rotatebox{180}{$\xrightarrow{\hspace{15pt}}$}}
\put(119,44){\scalebox{1.5}[1]{$\sim$}}
\put(113,56){\tiny $\varphi^{\alg{A}^-}_{F_1 \cap F_2}$}

\put(140,45){$\alg{A}^-\!/(F_1 \cap F_2)$}

\put(27,38){\rotatebox{270}{$\xrightarrow{\hspace{15pt}}$}}
\put(24,35){\rotatebox{270}{\scalebox{1.5}[1]{$\sim$}}}
\put(33,27){\tiny $\varphi^{\alg{B}^-}_{F_1 \cap B}$}

\put(97,38){\rotatebox{270}{$\xrightarrow{\hspace{15pt}}$}}
\put(97,39){\rotatebox{270}{$\xrightarrow{\hspace{14pt}}$}}
\put(103,28){\tiny $(i_Y)^*$}

\put(163,38){\rotatebox{270}{$\xrightarrow{\hspace{15pt}}$}}
\put(163,39){\rotatebox{270}{$\xrightarrow{\hspace{14pt}}$}}
\put(169,28){\tiny $q$}

\put(5,85){$(\alg{B}/(\theta|_B))^-$}

\put(64,85){$\lhook$}
\put(65,85){$\xrightarrow{\hspace{75pt}}$}
\put(97,93){\tiny $j$}

\put(150,85){$(\alg{A}/\theta)^-$}

\put(27,78){\rotatebox{270}{$\xrightarrow{\hspace{15pt}}$}}
\put(24,75){\rotatebox{270}{\scalebox{1.5}[1]{$\sim$}}}
\put(33,65){\tiny $i_2$}

\put(163,78){\rotatebox{270}{$\xrightarrow{\hspace{15pt}}$}}
\put(160,75){\rotatebox{270}{\scalebox{1.5}[1]{$\sim$}}}
\put(169,65){\tiny $i_1$}

\end{picture}
\end{center}
\end{claim}

\begin{proof}
The map
$j\colon \alg{B}/(\theta|_B) \mrig \alg{A}/\theta$,
defined by $b/(\theta|_B) \mapsto b/\theta$, is an injective $\cl{K}$-epimorphism, by Lemma~\ref{lem:epi}.
By Lemma~\ref{lem:dualityProperties}(\ref{lem:dualityProperties:minus}), the restriction of $j$ to $(B/(\theta|_B))^-$ is a
$\cl{BRA}$-morphism from $(\sbB/(\theta|_B))^-$
into $(\alg{A}/\theta)^-$.  We shall not distinguish notationally between $j$ and this restriction.
Whenever $b\in B$ and $b/(\theta|_B)\leqslant e/(\theta|_B)$, then $b/(\theta|_B)=(b\wedge e)/(\theta|_B)$ and
$b/\theta=(b\wedge e)/\theta$, so
\begin{equation}\label{b-}
(B/(\theta|_B))^-=\{b/(\theta|_B):b\in B^-\}
\end{equation}
and $j[(B/(\theta|_B))^-]=\{b/\theta:b\in B^-\}$.

Let $K=F_1\cap F_2$, so $\theta=\leibniz^\sbA\Fg^\sbA K$.
By (\ref{going up}), $A^-\cap\Fg^{\alg{A}}K=K$, so $B^-\cap\Fg^\sbA K=B^-\cap K=B\cap K=B\cap F_1$
(since $B\cap F_1= B\cap F_2$).
By (\ref{omega and subalgebras}), $\theta|_B = \leibniz^{\alg{B}}(B\cap \Fg^{\alg{A}}K)$, so
Lemma~\ref{lem:dualityProperties}(\ref{lem:dualityProperties:FilterRestriction}) supplies isomorphisms
\[
i_1 \colon (\alg{A}/\theta)^-
\cong \alg{A}^-/(F_1 \cap F_2) \ \text{ and } \ i_2 \colon (\alg{B}/(\theta|_B))^-
\cong \alg{B}^-/(F_1 \cap B),
\]
defined by $a/\theta
\mapsto (a \wedge e)/(F_1 \cap F_2)$ and $b/(\theta|_B)
\mapsto (b \wedge e)/(F_1 \cap B)$.

By Lemma~\ref{lem:subspace}, $\upclose^{\alg{B}^-}\!(F_1 \cap B)$ is the universe of an E-subspace, $\spc{Z}$ say, of ${\alg{B}^-}_{\!*}$, and
$\varphi^{\alg{B}^-}_{F_1 \cap B}\colon
\alg{B}^-\!/(F_1 \cap B) \cong \spc{Z}^*$.  Also,
$q\colon a/(F_1 \cap F_2) \mapsto a/F_1$ defines a homomorphism from $\alg{A}^-\!/(F_1 \cap F_2)$ onto $\alg{A}^-\!/F_1$.
Let $\spc{X}$ [resp.\ $\spc{Y}$] be the E-subspace of ${\alg{A}^-}_{\!*}$ with universe $\upclose^{\alg{A}^-}\!(F_1 \cap F_2)$ [resp.\ $\upclose^{\sbA^-} \!F_1$].
Let $i_Y\colon
\spc{Y}\mrig
\spc{X}$ be the inclusion morphism in $\cl{PESP}$. By Lemma~\ref{lem:subspace}, the following diagram commutes.
\begin{center}
\begin{picture}(120,65)


\put(15,5){$\alg{A}^-\!/F_1$}

\put(50,6){$\xrightarrow{\hspace{40pt}}$}
\put(56,3){\scalebox{4}[1]{$\sim$}}
\put(67,15){\tiny $\varphi^{\alg{A}^-}_{F_1}$}

\put(101,5){$\spc{Y}^*$}

\put(30,38){\rotatebox{270}{$\xrightarrow{\hspace{15pt}}$}}
\put(30,39){\rotatebox{270}{$\xrightarrow{\hspace{14pt}}$}}
\put(37,28){\tiny $q$}

\put(103,38){\rotatebox{270}{$\xrightarrow{\hspace{15pt}}$}}
\put(103,39){\rotatebox{270}{$\xrightarrow{\hspace{14pt}}$}}
\put(110,28){\tiny $(i_Y)^*$}

\put(0,45){$\alg{A}^-\!/(F_1 \cap F_2)$}

\put(67,44){$\xrightarrow{\hspace{25pt}}$}
\put(73,41){\scalebox{2}[1]{$\sim$}}
\put(68,53){\tiny $\varphi^{\alg{A}^-}_{F_1 \cap F_2}$}

\put(101,45){$\spc{X}^*$}

\end{picture}
\end{center}

Recall that the $\cl{BRA}$-morphism $i\colon\sbB^-\mrig\sbA^-$ is the inclusion map.
As $i$ is injective, its dual $i_* \colon {\alg{A}^-}_{\!*} \mrig {\alg{B}^-}_{\!*}$ is surjective, by Lemma~\ref{lem:BRAProperties}(\ref{lem:BRAProperties:SurjectiveVSInjective}).

The above definitions clearly imply that $i_*[Y]\subseteq Z$.  To establish the reverse inclusion,
let $G \in Z$.  By the Prime Filter Extension Theorem, $G = H \cap B$ for some $H \in \Pr(\alg{A}^-)$. Now,
$i_*[F_1]=B\cap F_1 \subseteq
G=i_*[H]$, so by (\ref{esakia morphism}),
$i_*[H] = i_*[H']$ for some $H'\in\upclose^{\sbA^-}\!F_1=Y$, whence
$G=i_*[H']\in i_*[Y]$.  Therefore, $Z = i_*[Y]$.

We claim that $i_*|_Y$ is injective.
Suppose, on the contrary, that ${H_1,H_2 \in Y}$, with $H_1 \neq H_2$ and $i_*[H_1] = i_*[H_2]$.
For each $k\in \{1,2\}$, (\ref{eq:DefOfF1}) shows that $\depth(F_1) \leq \depth(H_k)$, but
$F_1 \subseteq
H_k$, so $H_k=F_1$, whence $H_1=H_2$.  This contradiction confirms that $i_*|_Y$ is injective,
whence $i_*|_Y\colon\spc{Y}\cong\spc{Z}$ in $\cl{PESP}$.  In $\cl{BRA}$, therefore,
$({i_*|_Y})^*\colon\spc{Z}^*\cong\spc{Y}^*$.

A composition of isomorphisms in $\cl{BRA}$ is an isomorphism, so
\begin{equation}\label{g}
g\seteq (i_*|_Y)^* \circ \varphi^{\alg{B}^-}_{F_1 \cap B}\,\colon\alg{B}^-\!/(F_1 \cap B)\, \cong \,\spc{Y}^*.
\end{equation}
To show that the diagram in Claim~\ref{claim:retract} commutes, it remains to prove that
$g\circ i_2=\varphi^{\sbA^-}_{F_1}\!\circ q\circ i_1\circ j$.
And indeed, if $b\in B$ and $b/(\theta|_B)\in (B/(\theta|_B))^-$, then
\begin{align*}
& (g\circ i_2)(b/(\theta|_B))=\,g((b\wedge e)/(F_1 \cap B)) \\
& \quad\quad\quad\quad\quad\quad = \,(i_*|_Y)^*(\{ H \in \Pr(\alg{B}^-) : (F_1\cap B)\cup\{b\wedge e\}\subseteq H\})\\
& \quad\quad\quad\quad\quad\quad = \,\{ H \in \Pr(\alg{A}^-) : F_1\subseteq H \textup{ and }(F_1\cap B)\cup\{b\wedge e\}\subseteq H\cap B
\}\\
& \quad\quad\quad\quad\quad\quad = \,\{ H \in \Pr(\alg{A}^-) : F_1\cup\{b\wedge e\}\subseteq H
\}\\
& \quad\quad\quad\quad\quad\quad = \,\varphi^{\alg{A}^-}_{F_1}\!((b \wedge e)/F_1)\,
= \,(\varphi^{\sbA^-}_{F_1}\!\circ q\circ i_1\circ j)(b/(\theta|_B)).\qedhere
\end{align*}

\end{proof}

\smallskip

\begin{claim}\label{claim:another bloody claim}
\emph{Suppose\/ $k\in\{1,2\}$ and\/ $a\in A^-$ and\/ $b\in B^-$\textup{,} where\/
$a\equiv_\theta b
$\textup{.}
Then\/ $a\in F_k$ iff\/ $b\in F_k$\textup{.}  Consequently,
$a\notin (F_1\ld F_2)\cup(F_2\ld F_1)$\textup{.}}
\end{claim}
\begin{proof}
As $
a\equiv_\theta b
$, we have $a\leftrightarrow b\in\Fg^\sbA(F_1\cap F_2)$, so
\[
a\leftrightarrow^-b\seteq (a\leftrightarrow b)\wedge e\in A^-\cap\Fg^\sbA(F_1\cap F_2)=F_1\cap F_2,
\]
by (\ref{going up}).
As $a\rig^-b,\,b\rig^-a\geqslant a\leftrightarrow^- b$, it follows that $a\rig^-b,\,b\rig^- a\in F_k$.
Thus, $a\in F_k$ iff $b\in F_k$.  In particular, if $a\in (F_1\ld F_2)\cup(F_2\ld F_1)$ then
$b\in B\cap ((F_1\ld F_2)\cup(F_2\ld F_1))$, contradicting the fact that $B\cap F_1=B\cap F_2$.
\end{proof}

By Claim~\ref{claim:retract},
$h \seteq \varphi^{\alg{A}^-}_{F_1 \cap F_2} \circ i_1\,\colon (\alg{A}/\theta)^-\cong\spc{X}^*$ and, for each $a\in A$ such that
$a/\theta\in (A/\theta)^-$, we have
\[
h(a/\theta)=
\{ H \in \Pr(\alg{A}^-) :
(F_1 \cap F_2)\cup\{a\wedge e\} \subseteq H \}.
\]

By Claim~\ref{claim:2possibilities}, $F_1\ld F_2\neq \emptyset$.  In fact,
\begin{equation*}
\textup{$h(d/\theta) = \upclose^{\sbA^-}\!F_1$ for all $d\in F_1\ld F_2$.}
\end{equation*}
To confirm this, let $d\in F_1\ld F_2$.
Clearly, $\upclose^{\sbA^-}\!F_1\subseteq h(d/\theta)$.
Conversely, let $H \in h(d/\theta)$.  If $F_1\not\subseteq H$, then
since $F_1 \cap F_2 \subseteq H\in\Pr(\sbA^-)$,
we have
$F_2 \subseteq H$, by
(\ref{lem:filterProperties:ForG}).
In that case,
$F_2 \subsetneq
H$, because $d \in H \ld F_2$, but this contradicts
Lemma~\ref{lem:F2ltG}.  Thus,
$F_1 \subseteq H$, and so $h(d/\theta) = \upclose^{\sbA^-}\!F_1$, as claimed.

In Case (\ref{claim:2possibilities:F2ltF1}) of Claim~\ref{claim:2possibilities},
we have
$\upclose^{\sbA^-}\!F_2 = \upclose^{\sbA^-}\!(F_1\cap F_2)=X = h(e/\theta)$.

In Case (\ref{claim:2possibilities:incomparable}), we have $F_2\ld F_1\neq\emptyset$,
and we claim that
\begin{equation*}
h(d/\theta
) = \upclose^{\sbA^-}\!F_2 \textup{ \,for every $d\in F_2\ld F_1$.}
\end{equation*}
To see this, let $d\in F_2\ld F_1$.
It is clear that $\upclose^{\sbA^-}\!F_2 \subseteq h(d/\theta
)$, so consider $H \in h(d/\theta
)$.  If $F_2\not\subseteq H$ then, since $F_1\cap F_2\subseteq H$, we have
$F_1 \subseteq H$, by (\ref{lem:filterProperties:ForG}).  In that case, $F_1 \subsetneq
H$, as $d \in H \ld F_1$, but this contradicts Lemma~\ref{lem:F1ltG}.
Thus, $F_2\subseteq H$, and so $h(d/\theta
) \subseteq \upclose^{\sbA^-}\!F_2$.

We define $M=(A/\theta)^- \setminus j[(B/(\theta|_B))^-]$.

\smallskip

\begin{claim}
\label{claim:cardinality}
\emph{Fix any\/ $a_1 \in F_1 \ld F_2$\textup{.}  Choose\/ $a_2$ to be\/ $e$ in Case~\textup{(\ref{claim:2possibilities:F2ltF1})}
of Claim~\textup{\ref{claim:2possibilities},} and an arbitrary element of\/ $F_2\ld F_1$ in Case~\textup{(\ref{claim:2possibilities:incomparable}).}
Then\/
\[
\text{$M=\{a_1/\theta\}$ \,in Case~\textup{(\ref{claim:2possibilities:F2ltF1}),} \,and \,$M=\{a_1/\theta,\,a_2/\theta\}$
\,in Case~\textup{(\ref{claim:2possibilities:incomparable}).}}
\]
Moreover,
$a/\theta \prec e/\theta$ (in\/ $\sbA/\theta$) for all\/ $a\in A$ such that\/ $a/\theta\in M$\textup{.}}
\end{claim}

\begin{proof}
Observe that $a_1,a_2\leqslant e$ and, as we showed above,
\[
\textup{$h(a_1/\theta) = \upclose^{\sbA^-}\!F_1$ \,and\,
$h(a_2/\theta) = \upclose^{\sbA^-}\!F_2$.}
\]
By Claim~\ref{claim:another bloody claim} and (\ref{b-}), we have
$a_1/\theta\in M$ and, in Case~(\ref{claim:2possibilities:incomparable}), $a_2/\theta\in M$.
In Case~(\ref{claim:2possibilities:F2ltF1}),
$a_2/\theta\notin M$, since $e\in B$.

Because $h$ is an isomorphism,
$h[M]= X^*
\setminus h[j[(B/(\theta|_B))^-]]$ and,
for the first assertion of Claim~\ref{claim:cardinality}, it suffices to prove that
$h[M]\subseteq \{ \upclose^{\sbA^-}\!F_1,\, \upclose^{\sbA^-}\!F_2 \}$.

Suppose, with a view to contradiction, that there exists
$\mathcal{U} \in h[M]$ with
$\mathcal{U} \notin \{\upclose^{\sbA^-}\!F_1,\, \upclose^{\sbA^-}\!F_2\}$. Then
$\mathcal{U}\subseteq X$, but $\mathcal{U} \neq X$, because
\[
X = h(j(e/(\theta|_B)))\in h[j[(B/(\theta|_B))^-]].
\]
We show first that $\mathcal{U} \subsetneq \upclose^{\sbA^-}\!F_1$.

Suppose $F_1,F_2\in \mathcal{U}$.  For each $H \in X$, we have $F_1 \subseteq H$ or $F_2 \subseteq H$, by (\ref{lem:filterProperties:ForG}), so
$H \in \mathcal{U}$ (since $\mathcal{U}$ is upward closed).
This shows that $X\subseteq\mathcal{U}$, a contradiction. Therefore, $F_1$ and $F_2$ don't both belong to $\mathcal{U}$.

Suppose $F_2 \in \mathcal{U}$. Then $F_1 \notin \mathcal{U}$ and $\upclose^{\sbA^-}\!F_2 \subseteq \mathcal{U}$, as $\mathcal{U}$ is upward closed.  If $H \in \mathcal{U}$, then $H\neq F_1$ and $F_1 \cap F_2 \subseteq H$ (as $\mathcal{U} \subseteq X$). In that case, $F_2 \subseteq H$ (otherwise, $F_1 \subseteq H$, by
(\ref{lem:filterProperties:ForG}), whence $F_1 \subsetneq
H$, but then $F_2 \subsetneq
H$, by Lemma~\ref{lem:F1ltG}).  This shows that $\mathcal{U} \subseteq \upclose^{\sbA^-}\!F_2$, so $\mathcal{U} = \upclose^{\sbA^-}\!F_2$,
contrary to our initial assumptions about $\mathcal{U}$.
Therefore, $F_2 \notin \mathcal{U}$.

We claim that $\mathcal{U} \subseteq \upclose^{\sbA^-}\!F_1$.  For otherwise, $F_1 \not\subseteq H$ for some $H \in \mathcal{U}$, whence $H\neq F_2$
and, by (\ref{lem:filterProperties:ForG}), $F_2 \subseteq H$, i.e., $F_2 \subsetneq
H$, whereupon Lemma~\ref{lem:F2ltG} delivers the contradiction
$F_1 \subseteq H$.
Thus, $\mathcal{U} \subsetneq
\upclose^{\sbA^-}\!F_1$ (since $\mathcal{U} \neq \upclose^{\sbA^-}\!F_1$, by assumption).

Now we shall argue that $\mathcal{U} \in  h[j[(B/\theta|_B)^-]]$ (contradicting the fact that
$\mathcal{U}\in h[M]$, and thereby  confirming the relation
$h[M]\subseteq \{ \upclose^{\sbA^-}\!F_1,\, \upclose^{\sbA^-}\!F_2 \}$).

As $\mathcal{U}\in X^*$ and $\mathcal{U} \subsetneq
\upclose^{\sbA^-}\!F_1 = Y$, we have
\[
\mathcal{U} = \mathcal{U} \cap Y = (i_Y)^*(\mathcal{U})\in Y^*,
\]
so by Claim~\ref{claim:retract}, there exists $b\in B$ with $b/(\theta|_B)\in (B/(\theta|_B))^-$ such that
\[
\mathcal{U} = g(i_2(b/(\theta|_B))) = (i_Y)^*(h(j(b/(\theta|_B)))) = Y\cap h(b/\theta)=Y\cap\mathcal{V},
\]
where $g$ is
as in (\ref{g}), and
$\mathcal{V} \seteq h(b/\theta)$.
By (\ref{b-}), we may assume
that $b\leqslant e$.

Now $\mathcal{V}\in X^*$, so $\mathcal{V}$ is an upward-closed subset of $X$.
Note that $F_1 \notin \mathcal{V}$ (otherwise
$Y=\upclose^{\sbA^-}\!F_1 \subseteq \mathcal{V}$, yielding the contradiction
$\mathcal{U}\subsetneq Y=Y\cap\mathcal{V}=\mathcal{U}$).  It follows that $Y\not\subseteq\mathcal{V}$ (as $F_1\in Y$).

For any $H \in \mathcal{V}$, if $F_1 \subseteq H$, then $F_1 \subsetneq
H$ (as $F_1\notin\mathcal{V}$), whence $F_2 \subsetneq
H$ (by Lemma~\ref{lem:F1ltG}), whereas if $F_1\not\subseteq H$, then $F_2\subseteq H$ (by (\ref{lem:filterProperties:ForG})).
This shows that
$\mathcal{V} \subseteq \upclose^{\sbA^-}\!F_2$.

We now argue that $\mathcal{V}\subseteq Y$.

Suppose, on the contrary, that there exists $H \in \mathcal{V} \ld Y$. Then $F_1 \not\subseteq H$ (by definition of $Y$), so $F_2 \subseteq H$,
by (\ref{lem:filterProperties:ForG}).  Now Lemma~\ref{lem:F2ltG} prevents $F_2$ from being a proper subset of $H$, so
$F_2 = H$.
In particular, $F_2\in\mathcal{V}$, so $\upclose^{\sbA^-}\!F_2 \subseteq \mathcal{V}$, whence
$\mathcal{V} = \upclose^{\sbA^-}\!F_2$.

In Case (\ref{claim:2possibilities:F2ltF1}) of Claim~\ref{claim:2possibilities}, it would follow that
$Y=\upclose^{\sbA^-}\!F_1\subseteq \upclose^{\sbA^-}\!F_2=\mathcal{V}$, a contradiction.

In Case~(\ref{claim:2possibilities:incomparable}), we have $e\geqslant a_2 \in F_2 \ld F_1$ and
$h(a_2/\theta) = \upclose^{\sbA^-}\!F_2=\mathcal{V}=h(b/\theta)$.
Then, since $h$ is injective, $a_2/\theta = b/\theta$, contradicting Claim~\ref{claim:another bloody claim}.

This confirms that
$\mathcal{V} \subseteq Y$, and so $\mathcal{V}=Y\cap \mathcal{V}
= \mathcal{U}$.  Therefore,
\[
\mathcal{U} = h(b/\theta)=h(j(b/(\theta|_B)))\in h[j[(B/\theta|_B)^-]],
\]
completing the proof that
$M$ is $\{a_1/\theta\}$ in Case~(\ref{claim:2possibilities:F2ltF1}), and is $\{a_1/\theta,\,a_2/\theta\}$
in Case~(\ref{claim:2possibilities:incomparable}).

\smallskip

It remains to show that $a/\theta \prec e/\theta$ in $\sbA/\theta$, whenever $a/\theta\in M$.

To establish that $a_1/\theta\prec e/\theta$ in $\sbA/\theta$ (i.e., in $(\sbA/\theta)^-$), it suffices
to show that $\upclose^{\sbA^-}\!F_1\prec X$ in $\sbX^*$, because $h$ is an isomorphism.

Suppose
$\upclose^{\sbA^-}\!F_1 \subsetneq
\mathcal{W} \subsetneq
X$, where $\mathcal{W} \in X^*$.  Then
$F_1 \not\subseteq H$ for some
$H \in \mathcal{W}$, whence
$F_2 \subseteq H$, by (\ref{lem:filterProperties:ForG}).
We cannot have $F_2 = H$, otherwise $\upclose^{\sbA^-}\! F_2 \subseteq \mathcal{W}$, in which case every
element $G$ of $X$ belongs to $\mathcal{W}$ (as $G$ contains $F_1$ or $F_2$, again by (\ref{lem:filterProperties:ForG})).
Therefore,
$F_2 \subsetneq
H$, and so $F_1 \subseteq H$, by Lemma~\ref{lem:F2ltG}.  This contradiction confirms that
$a_1/\theta \prec e/\theta$ in $\sbA/\theta$.

We may now assume that Case~(\ref{claim:2possibilities:incomparable}) applies.  The desired conclusion $\textup{$a_2/\theta\prec e/\theta$}$
amounts similarly to the claim that $\upclose^{\sbA^-}\! F_2 \prec X$ in $\sbX^*$.  Suppose
$\upclose^{\sbA^-}\! F_2 \subsetneq \mathcal{W} \subsetneq X$, where
$\mathcal{W} \in X^*$.  Then $F_2\not\subseteq H$ for some $H\in\mathcal{W}$.  Now $H\in X$, so
$F_1\subseteq H$, by (\ref{lem:filterProperties:ForG}).  Then $F_1=H$, by Lemma~\ref{lem:F1ltG}, so $F_1\in\mathcal{W}$,
whence $\upclose^{\sbA^-}\!F_1\subseteq\mathcal{W}$.  By (\ref{lem:filterProperties:ForG}) again,
$X\subseteq (\upclose^{\sbA^-}\!F_1)\cup(\upclose^{\sbA^-}\!F_2)\subseteq \mathcal{W}$, contradicting
the fact that $\mathcal{W}\subsetneq X$.
Therefore, $a_2/\theta \prec e/\theta$ in $\sbA/\theta$.
\end{proof}

We are now in a position to prove Claim~\ref{claim:subalgebra} (and hence Theorem~\ref{thm:SIRLhasES}).

\begin{proof}[Proof of Claim~\ref{claim:subalgebra}]
Since $\sbA$ and $\sbB$ are negatively generated,
so are $\sbA/\theta$ and
$\sbB/(\theta|_B)$, by Lemma~\ref{neg cone gen in images}.  The subalgebra
\[
\sbJ\seteq j[\sbB/(\theta|_B)]
\]
of $\sbA/\theta$ is isomorphic to
$\sbB/(\theta|_B)$, so $\sbJ$ is also negatively generated.
By
Lemma~\ref{lem:dualityProperties}(\ref{lem:dualityProperties:minus}),
$J^-=j[(B/(\theta|_B))^-]$, whence $M=(A/\theta)^-\ld J^-$.
As $\sbB$ is $\cl{K}$-epic in $\sbA$,
Lemma~\ref{lem:epi} shows that $\sbJ$ is $\cl{K}$-epic in
$\sbA/\theta$.
Moreover,
\[
\sbS\seteq \boldsymbol{\mathit{Sg}}^{\sbA/\theta}(J^-\cup\{a_1/\theta\})
\]
is negatively generated
(since $J^-\cup\{a_1/\theta\}\subseteq S^-$),
and $\sbJ$ is a subalgebra of $\sbS$ (as $J=\Sg^{\sbA/\theta}(J^-)$), so
$\sbS$ is $\cl{K}$-epic in $\sbA/\theta$ (because $\sbJ$ is).

Observe that $J\neq A/\theta$, because $a_1/\theta\notin J$ (by Claim~\ref{claim:another bloody claim} and (\ref{b-}), since
$a_1\in F_1\ld F_2$),
and that $A/\theta=\Sg^{\sbA/\theta}((A/\theta)^-)=\Sg^{\sbA/\theta}(J^-\cup M)$.

We choose $\sbC=\sbJ$ and $a=a_1$ in Case~(\ref{claim:2possibilities:F2ltF1}).
We make the same choices in Case~(\ref{claim:2possibilities:incomparable}) \emph{if} $a_2/\theta\in S$.
Under these conditions, $J^-\cup M=J^-\cup\{a_1/\theta\}$ (by Claim~\ref{claim:cardinality}) and
$A/\theta=
\Sg^{\sbA/\theta}(J^-\cup M)\subseteq S=\Sg^{\sbA/\theta}(C^-\cup\{a/\theta\})$,
so $\sbA/\theta$ is generated by $C^-\cup\{a/\theta\}$, as required.

In Case~(\ref{claim:2possibilities:incomparable}), if $a_2/\theta\notin S$ (whence $S\neq A/\theta$), we choose $\sbC=\sbS$ and $a=a_2$, whereupon $J^-\cup M=J^-\cup\{a_1/\theta,a_2/\theta\}$ (by Claim~\ref{claim:cardinality}) and
\[
A/\theta=
\Sg^{\sbA/\theta}(J^-\cup M)\subseteq \Sg^{\sbA/\theta}(S^-\cup\{a_2/\theta\})=\Sg^{\sbA/\theta}(C^-\cup\{a/\theta\}),
\]
so again, $\sbA/\theta$ is generated by $C^-\cup\{a/\theta\}$.
\end{proof}

\section{Reflections and De Morgan Monoids}\label{reflections}

Given an SRL $\sbA$,
let $A'=\{a':a\in A\}$ be a disjoint copy of $A$, and let $\bot,\top$ be distinct non-elements of $A\cup A'$.
The \emph{reflection} $\textup{R}(\sbA)$ of $\sbA$ is the SIRL with universe $\textup{R}(A)=A\cup A'\cup\{\bot,\top\}$ such that
$\sbA$ is a subalgebra of the SRL-reduct of $\textup{R}(\sbA)$ and, for all $a,b\in A$ and $x,y\in \textup{R}(A)$,
\begin{align*}
& x\bcdw\bot=\bot < a < b' < \top = a'\bcdw b', \textup{ and if $x\neq\bot$, then }x\bcdw\top=\top;\\
& a\bcdw b' = (a\rig b)';\\
& \neg a = a' \textup{ and } \neg(a')=a \textup{ and } \neg\bot=\top \textup{ and } \neg\top = \bot;\\
& x\rig y=\neg(x\bcdw\neg y).
\end{align*}
Since $f=e'$, we have $\top=f^2$ and $\bot=\neg(f^2)$, so $\bot,\top$ belong to every subalgebra of $\textup{R}(\sbA)$.

The reflection construction originates with Meyer \cite{Mey73}; also see \cite[Section~9]{GR04}.  It preserves (and reflects) distributivity and the square-increasing law, so
$\sbA$ is a Dunn monoid iff $\textup{R}(\sbA)$ is a De Morgan monoid.
Also, $\sbA$ is FSI iff $\textup{R}(\sbA)$ is.

The \emph{reflection}
of a variety $\mathsf{K}$ of SRLs
is the
variety
\[
\mathbb{R}(\mathsf{K})\seteq\mathbb{V}\{\textup{R}(\sbA):\sbA\in\mathsf{K}\}
\]
of SIRLs.
The structure of a member of $\mathbb{R}(\cl{K})$ is illuminated by the following lemma, which is proved in \cite[Section~6]{MRW2}.
(The extra
assumptions there were not relied on in the proof.)

\begin{lemma}\label{hs lem}
\textup{(cf.\ \cite[Lemma~6.5]{MRW2})}\,
Let\/ $\sbA$ be an SRL.
\begin{enumerate}
\item\label{s}
If\/ $\sbB$ is a subalgebra of $\sbA$\textup{,} then\/
$B\cup\{b':
b\in B\}\cup\{\bot,\top\}$
is the universe of a subalgebra of\/ $\textup{R}(\sbA)$ that is isomorphic to $\textup{R}(\sbB)$\textup{,} and every
subalgebra of\/ $\textup{R}(\sbA)$ arises in this way from a subalgebra of\/ $\sbA$\textup{.}

\smallskip
\item\label{h}
If\/ $\theta$ is a congruence of\/ $\sbA$\textup{,} then
\[
\textup{\quad\quad $\textup{R}(\theta)\seteq\theta\cup\{\langle a',b'\rangle:
\langle a,b\rangle\in \theta\}\cup\{\langle \bot,\bot\rangle,\,\langle \top,\top\rangle\}$}
\]
is a congruence of\/ $\textup{R}(\sbA)$\textup{,} and\/ $\textup{R}(\sbA)/\textup{R}(\theta)\cong\,\textup{R}(\sbA/\theta)$\textup{.}  Also, every proper congruence of\/
$\textup{R}(\sbA)$ has the form\/ $\textup{R}(\theta)$ for some\/ $\theta\in{{\mathit{Con}}}\,\sbA$\textup{.}

\smallskip

\item\label{pu}
If\/ $\{\sbA_i:
i\in I\}$ is a family of SRLs
and\/ $\mathcal{U}$ is an ultrafilter over\/ $I$\textup{,} then\/
$\prod_{i\in I}\textup{R}(\sbA_i)/\mathcal{U}\,\cong\textup{R}\!\left(\prod_{i\in I}\sbA_i/\mathcal{U}\right)$\textup{.}
\end{enumerate}
\end{lemma}

\emph{J\'{o}nsson's Theorem} \cite{Jon67,Jon95} states that, for any subclass
$\mathsf{L}$ of a
congruence distributive variety,
$\mathbb{V}(\mathsf{L})_{FSI}
\subseteq
\mathbb{HSP}_\mathbb{U}(\mathsf{L})$.  Together with
Lemma~\ref{hs lem}, this yields the next corollary (as every variety is generated by
its FSI members).

\begin{corollary}\label{fsi in reflections}
Let\/ $\cl{K}$ be a variety of SRLs, with $\sbE\in\mathbb{R}(\cl{K})$\textup{.}  Then\/ $\sbE$ is FSI iff
$\sbE\cong\textup{R}(\sbD)$ for some\/ $\sbD\in\cl{K}_{FSI}$\textup{.}
\end{corollary}

\begin{theorem}\label{ES and reflections}
Let\/ $\cl{K}$ be a variety of SRLs, let\/ $\sbB$ be a subalgebra of\/ $\sbA\in\cl{K}$\textup{,}
and identify\/ $\textup{R}(\sbB)$ with the subalgebra of\/ $\textup{R}(\sbA)$ given in Lemma~\textup{\ref{hs lem}(\ref{s}).}
Then\/
\begin{enumerate}
\item\label{i}
$\sbB$ is\/ $\cl{K}$-epic in\/ $\sbA$ iff\/ $\textup{R}(\sbB)$ is\/ $\mathbb{R}(\cl{K})$-epic in\/ $\textup{R}(\sbA)$\textup{;}

\item\label{ii}
$\cl{K}$ has the ES property iff\/ $\mathbb{R}(\cl{K})$ has the ES property;

\item\label{iii}
$\cl{K}$ is locally finite iff\/ $\mathbb{R}(\cl{K})$ is locally finite.
\end{enumerate}
\end{theorem}
\begin{proof}
(\ref{i}) ($\Rig$): Let $g,h\colon\textup{R}(\sbA)\mrig\sbE\in\mathbb{R}(\cl{K})$ be homomorphisms that agree on $\textup{R}(\sbB)$.
In showing that $g=h$, we may assume that $\sbE$ is subdirectly irreducible (by the Subdirect Decomposition Theorem),
whence $\sbE=\textup{R}(\sbD)$ for some $\sbD\in\cl{K}_{FSI}$, by Corollary~\ref{fsi in reflections}.
Since $g,h$ preserve $e,\bcdw,\neg$, they preserve $\bot,\top$.  If $a,b\in A$, then $g(a),h(a)\neq\bot$
(otherwise, the kernel of $g$ or $h$ would identify $\top=a\bcdw\top$ with $\bot\bcdw\top=\bot$),
and $g(a),h(a)\neq\top$ (because the kernels don't identify $\top=\top\rig\top$ with $\top\rig a=\bot$),
while $g(a),h(a)\neq b'$ (because the kernels don't identify $a^2\in A$ with $\top=(b')^2$).  Thus,
$g[A],h[A]\subseteq D$, and so $g|_A,h|_A$ are homomorphisms from $\sbA$ to $\sbD$, which agree on $\sbB$.
As $\sbD\in\cl{K}$ and $\sbB$ is $\cl{K}$-epic in $\sbA$, we conclude that $g|_A=h|_A$.  Then $g|_{A'}=h|_{A'}$,
since $g,h$ preserve $\neg$.  Consequently, $g=h$.

($\Leftarrow$): Let $g,h\colon\sbA\mrig\sbD\in\cl{K}$ be homomorphisms that agree on $\sbB$.
Then $\textup{R}(\sbD)\in\mathbb{R}(\cl{K})$.
Let $\ov{g},\ov{h}\colon\textup{R}(\sbA)\mrig\textup{R}(\sbD)$ be the respective extensions of
$g,h$, preserving $\bot,\top$, such that $\ov{g}(a')=g(a)'$ and $\ov{h}(a')=h(a')$ for all $a\in A$.
Then $\ov{g},\ov{h}$ are homomorphisms that agree on $\textup{R}(\sbB)$, so by assumption, $\ov{g}=\ov{h}$,
whence $g=h$.

(\ref{ii}) Obviously, $\sbB=\sbA$ iff $\textup{R}(\sbB)=\textup{R}(\sbA)$.
Therefore, the implication from right to left follows from (\ref{i}).  For the converse,
use Theorem~\ref{thm:ESWitnessedByFSI}, Corollary~\ref{fsi in reflections},
Lemma~\ref{hs lem}(\ref{s}) and item~(\ref{i}) of the present theorem.

(\ref{iii}) ($\Rig$): As $\cl{K}$ is locally finite, there is a function $p\colon\omega\mrig\omega$ such that,
for each $n\in\omega$, every $n$-generated member of $\cl{K}_{FSI}$ has at most $p(n)$ elements.  It
suffices to show that, for each $n\in\omega$, every $n$-generated $\sbE\in\mathbb{R}(\cl{K})_{FSI}$
has at most $2+2p(n)$ elements.  By Corollary~\ref{fsi in reflections}, any such $\sbE$ may be assumed to be
$\textup{R}(\sbD)$ for some $\sbD\in\cl{K}_{FSI}$.  Let $G$ be an irredundant generating set for $\sbE$, with
$\left|G\right|\leq n$.  Then $\bot,\top\notin G$.  Let $H=(G\cap D)\cup\{\neg g:g\in G\cap D'\}$, so $\left|H\right|\leq n$ and
$\sbC\seteq \Sg^\sbD H$ has at most $p(n)$ elements.  By Lemma~\ref{hs lem}(\ref{s}), $\textup{R}(\sbC)$
may be identified with a subalgebra of $\sbE$, but then $G\subseteq\textup{R}(C)$, so $\textup{R}(\sbC)=\sbE$,
whence $\left|E\right|\leq 2+2p(n)$.

($\Leftarrow$):  Use the fact that an SIRL of the form $\textup{R}(\sbA)$ is generated by $A$.
\end{proof}

As a function from the lattice of varieties of SRLs into that of SIRLs, the operator $\mathbb{R}$
is obviously isotone.  Using
Corollary~\ref{fsi in reflections}, we can show that
$\mathbb{R}$ is also $\subseteq$-reflecting (and therefore injective); the proof is the same as that of
\cite[Lemma~6.7]{MRW2} (where again, the extra assumptions play no role).

Consequently, as all varieties of Brouwerian algebras
of finite depth have the ES property (by \cite[Theorem~5.4]{BMR17} or by Theorem~\ref{thm:SIRLhasES})
and since there are $2^{\aleph_0}$ such distinct varieties
(even of depth $3$) \cite{Kuz75}, the following can be inferred from Theorem~\ref{ES and reflections}(\ref{ii}),(\ref{iii}).

\begin{theorem}\label{continuum}
There are\/ $2^{\aleph_0}$ distinct locally finite varieties of De Morgan monoids with the ES property.
\end{theorem}

\section{Further Examples and Applications}\label{section:applications}

\subsection{The weak ES property}\

\smallskip

The \emph{weak ES property} for a variety $\cl{K}$ rules out non-surjective $\cl{K}$-epimor\-phisms $h\colon\sbA\mrig\sbB$
in all cases where $\sbB$ is generated by the union of $h[A]$ and a \emph{finite} set.  By \cite[Theorem~5.4]{MRW3}, it
is equivalent to the demand that no finitely generated member of $\cl{K}$ has a $\cl{K}$-epic proper subalgebra.  In varieties
of logic, it amounts to the so-called \emph{finite Beth property} for the corresponding deductive system
\cite[Theorem~3.14, Corollary~3.15]{BH06} (also see \cite[Theorem~7.9]{MRW3}).

Therefore, by an argument of Kreisel \cite{Kre60},
every variety of Brouwerian (or Heyting) algebras has the weak ES property.  It follows from a result of Campercholi
\cite[Corollary~6.5]{Cam18} that, in any finitely generated variety with a majority term (e.g., one generated by a finite
lattice-based algebra), the weak ES property entails the ES property.  This provides a different explanation of the
slightly earlier finding that all finitely generated varieties of Brouwerian (or Heyting) algebras have the ES property
\cite[Corollary~5.5, Theorem~7.2]{BMR17}.

Every variety with the weak ES property and the amalgamation property has (a strong form of) the ES property;
see \cite{Isb66,KMPT83,Rin72} and \cite[Section~2.5.3]{Hoo01}.  Consequently, in all varieties of Brouwerian (or Heyting) algebras,
amalgamability entails epimorphism-surjectivity, and the amalgamable varieties of these kinds have
been classified completely by Maksimova; see \cite{GM05,Mak77a,Mak77b}.

The situation is different for varieties of (possibly non-integral)
S[I]RLs, as they may lack the weak ES property, even when they are finitely generated (see Section~\ref{hypotheses} below).

\subsection{Hypotheses of the main theorem}\label{hypotheses}\

\smallskip

If $\cl{K}$ and $\cl{L}$ are subvarieties of a congruence distributive
variety, with ${\cl{M}=\mathbb{V}(\cl{K}\cup\cl{L})}$,
then
$\cl{M}_{FSI}=\cl{K}_{FSI}\cup\cl{L}_{FSI}$, by J\'{o}nsson's Theorem.
The hypotheses in Theorem~\ref{thm:SIRLhasES} therefore
persist in (binary) varietal joins, and of course in subvarieties.  This is helpful in applications, because
the ES property itself is not hereditary.

Neither of the two hypotheses in Theorem~\ref{thm:SIRLhasES} can be dropped.

To see that the finite depth assumption cannot be dropped for varieties of SRLs, it suffices to exhibit a variety
of Brouwerian algebras (of infinite depth) without the ES property.  This was done in \cite[Section~6]{BMR17}.
It was subsequently shown in \cite{MW} that there are $2^{\aleph_0}$ distinct locally finite varieties of Brouwerian
algebras that lack the ES property and that have \emph{width} $2$ (i.e., $2$ is the maximum cardinality of an anti-chain in
the dual of an FSI member of the variety).  It follows, as in Section~\ref{reflections}, that there are $2^{\aleph_0}$
locally finite varieties of
De Morgan monoids without the ES property, but this could alternatively be deduced from older findings, discussed below.

The demand for negative
generation is not redundant either, because some finitely generated
varieties of De Morgan monoids (and
of Dunn monoids) lack the ES property.
Indeed, by an argument of Urquhart \cite{Urq99} (also see \cite[Corollary~4.15]{BH06}), a variety of
De Morgan monoids lacks even the weak ES property if it
contains a certain six-element algebra $\sbC$, called the \emph{crystal lattice} (which is
not negatively generated).
That algebra is depicted below.  (Deletion of $b$ leaves
an epic subalgebra behind, owing to the uniqueness of existent relative complements
in distributive lattices.)  The argument adapts to Dunn monoids, using the SRL-reduct of $\sbC$.

{\tiny

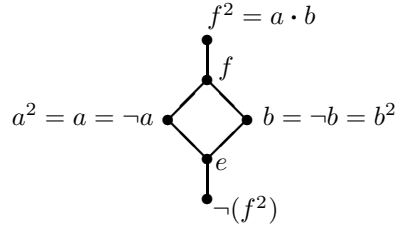
\thicklines
\begin{center}
\begin{picture}(80,90)(-130,39)

\put(-105,50){\circle*{4}}
\put(-105,65){\line(0,-1){15}}
\put(-105,65){\circle*{4}}
\put(-90,80){\line(-1,-1){15}}
\put(-90,80){\circle*{4}}
\put(-120,80){\line(1,-1){15}}
\put(-120,80){\circle*{4}}
\put(-120,80){\line(1,1){15}}
\put(-105,95){\circle*{4}}
\put(-90,80){\line(-1,1){15}}
\put(-105,95){\line(0,1){15}}
\put(-105,110){\circle*{4}}

\put(-105,116){\small $f^2=a\bcdw b$}
\put(-101,97){\small $f$}
\put(-179,78){\small $a^2=a = \neg a$}
\put(-84,78){\small $b = \neg b=b^2$}
\put(-102,61){\small $e$}
\put(-103,43){\small $\neg (f^2)$}

\end{picture}\nopagebreak
\end{center}
}

\noindent
In particular, $\sbC$ is absent from each of the $2^{\aleph_0}$ varieties $\cl{K}$ of De Morgan monoids
in Theorem~\ref{continuum}, while the corresponding locally finite
varieties $\mathbb{V}(\cl{K}\cup\{\sbC\})$ lack the weak ES property, and by J\'{o}nsson's Theorem,
they are distinct.

Although the finite depth assumption in Theorem~\ref{thm:SIRLhasES} cannot be dropped, it is not a \emph{necessary}
condition for the ES property.  Indeed, epimorphisms are surjective in $\cl{BRA}$ (hence in its reflection), and in the
locally finite variety of \emph{relative Stone algebras} (i.e., subdirect products of totally ordered Brouwerian algebras).
Also, the smallest variety containing the De Morgan monoids (alternatively,
the idempotent SRLs) that are totally ordered and negatively generated
has the ES property
and is locally finite \cite{Wan}.  All of the varieties mentioned in this paragraph have infinite depth.

Likewise, the demand that FSI members be negatively generated is not entailed by the ES property, even in
varieties of S[I]RLs of finite depth.  Witnessing examples can also be found in \cite{Wan}.

\subsection{More varieties of De Morgan monoids}\label{new positive}\

\smallskip

\emph{Sugihara monoids}---i.e., idempotent De Morgan monoids---are always negatively generated,
and the
same applies to the SRLs that can be
embedded into them (these are the \emph{positive Sugihara monoids} of \cite{OR07}).
The ES property was established recently for \emph{all} varieties of [positive] Sugihara monoids
\cite[Theorems~8.5,\,8.6]{BMR17}.
For the finitely generated varieties of this kind, the surjectivity of epimorphisms could
alternatively be deduced (immediately) from Theorem~\ref{thm:SIRLhasES}.

Apart from reflections and idempotent cases, Theorem~\ref{thm:SIRLhasES} yields further examples
as follows.

By \cite[Theorem~6.1]{MRW}, the lattice of varieties of De Morgan monoids has just four atoms, each of which is a finitely generated
variety satisfying the hypotheses of Theorem~\ref{thm:SIRLhasES}.  One of them is
generated by the reflection of a trivial SRL, i.e., by the non-idempotent De Morgan monoid $\sbC_4$ on the chain
\[
\neg(f^2)<e<f<f^2.
\]
The covers of $\mathbb{V}(\sbC_4)$
are distinctive, as $\sbC_4$
is the only $0$-generated nontrivial algebra
onto which FSI De Morgan monoids may be mapped by
non-injective homomorphisms \cite[Theorem~1]{Sla89}.

There is a largest variety $\mathsf{U}$
of De Morgan monoids consisting of homomorphic \emph{pre-images} of $\sbC_4$ (along with trivial algebras), and in
the subvariety lattice of $\mathsf{U}$, the variety $\mathbb{V}(\sbC_4)$ has just
ten covers \cite[Sections~4,\,8]{MRW2}, only two of which are generated by reflections of Dunn monoids.

Each of these ten varieties is generated by a finite De Morgan monoid
that is itself generated by one of the lower bounds of its neutral element (and is thus negatively generated).
Therefore, the conditions of Theorem~\ref{thm:SIRLhasES} obtain in all
ten covers, and hence in their varietal join, so all subvarieties of this join
have the ES property.

\subsection{Bounds}\label{bounds}\

\smallskip

The original Esakia duality of \cite{Esa85} supplies an equivalence between the category $\cl{HA}$ of
Heyting algebras (and their homomorphisms---which must preserve $\bot$) and the opposite of the category $\cl{ESP}$
of \emph{Esakia spaces}.

The objects of $\cl{ESP}$ are like those of $\cl{PESP}$, except that they need not have
maximum elements; the definition of morphisms is unaffected.  For $\sbA\in\cl{HA}$ and $\sbX\in\cl{ESP}$, we \emph{re-define}
$\Pr(\sbA)$ as the set of prime \emph{proper} filters of $\sbA$, and $\Cpu(\sbX)$ as the set of \emph{all} clopen up-sets
of $\sbX$, including $\emptyset$.  After these changes, the definitions of $\sbA_*$, $\sbX^*$, the duals of morphisms, and
the canonical isomorphisms remain the same (but note that $\sbA_*$ is empty when $\left|A\right|=1$).  The definition of depth
is adjusted so that a Heyting algebra and its Brouwerian reduct have the same depth.  (In particular, the depth of the Esakia
space reduct of a pointed
Esakia space $\sbX$ exceeds that of $\sbX$ by $1$.)

Theorem~\ref{thm:SIRLhasES} remains true for varieties of bounded S[I]RLs; its
proof requires no further alteration.  In this
form, it generalizes the
recent
finding
that every variety of Heyting
algebras of finite depth has surjective epimorphisms \cite[Theorem~5.3]{BMR17}.

\end{document}